\documentclass{article}

\usepackage{arxiv}

\usepackage[utf8]{inputenc} 
\usepackage[T1]{fontenc}    
\usepackage{hyperref}       
\usepackage{url}            
\usepackage{booktabs}       
\usepackage{amsfonts}       
\usepackage{nicefrac}       
\usepackage{microtype}      
\usepackage{algorithm2e}
\SetKwInput{KwData}{Input}
\SetKwInput{KwResult}{Output}

\usepackage{braket}
\usepackage{graphicx}
\usepackage{caption, subcaption}
\usepackage{amssymb,amsmath,color,soulutf8,longtable,colortbl,setspace,ifthen,xspace,url,pdflscape}
\usepackage{amsthm}
\usepackage{tikz}
\usetikzlibrary{3d}
\usetikzlibrary{shapes,arrows,positioning}
\usetikzlibrary{fit}
\usetikzlibrary{arrows.meta}
\usepackage{multirow}
\usepackage{comment}

\newcommand{\II}[0]{\mathbf{I}\hspace{-0.8pt}\mathbf{I}}
\newcommand{\bfu}[0]{\mathbf{u}}
\newcommand{\bfv}[0]{\mathbf{v}}

\title{Automatic feature-preserving size field\\for 3D mesh generation}

\author{
  Arthur Bawin,\\
  Institute of Mechanics, Materials and Civil Engineering\\
  Université catholique de Louvain\\
  and\\
  Département de génie mécanique\\
  École Polytechnique de Montréal\\
  \texttt{arthur.bawin@uclouvain.be} \\
  \And
 Fran{\c c}ois Henrotte \\
 Institute of Mechanics, Materials and Civil Engineering\\
 Université catholique de Louvain \\
 and\\
  ACE - Université de Liège\\
 \texttt{francois.henrotte@uclouvain.be} \\
   \And
 Jean-Fran{\c c}ois Remacle \\
  Institute of Mechanics, Materials and Civil Engineering\\
  Université catholique de Louvain \\
  \texttt{jean-francois.remacle@uclouvain.be} \\
}

\begin{document}
\maketitle

\begin{abstract}
This paper presents a methodology aiming at easing considerably
the generation of high-quality meshes for complex 3D domains.
We show that the whole mesh generation process
can be controlled with only five parameters
to generate in one stroke quality meshes for arbitrary geometries.
The main idea is to build a meshsize field $h({\bf x})$
taking local features of the geometry, such as curvatures, 
into account.
Meshsize information is then
propagated from the surfaces into the volume,
ensuring that the magnitude of $|\nabla h|$ is always controlled
so as to obtain a smoothly graded mesh.
As the meshsize field is stored in an independent octree data structure,
the function $h$ can be computed separately,
and then plugged in into any mesh generator
able to respect a prescribed meshsize field.
The whole procedure is automatic,
in the sense that minimal interaction with the user is required.
Applications examples based on models taken from the very large ABC dataset,
are then presented, all treated with the same generic set of parameter values,
to demonstrate the efficiency
and the universality of the technique.
\end{abstract}


\section{Introduction}
Models used in industry have considerably grown in complexity
over the last decades,
and it is now common to mesh models with tens of thousands of faces.
Ideally, 
a designer should create the CAD model,
press the \emph{generate mesh} button,
and obtain in less than a minute a computational mesh 
valid \emph{as is} for a 
finite element simulation.
Practitioners in the field know however
that things do not work out that easily in reality.
Mesh generation for complex geometries is in practice
a time-consuming task 
often involving intermediary meshes,
progressively enhanced to fulfill specified meshsize and quality
requirements.

The purpose of mesh generation is to build meshes
with elements of controlled size and quality.
We call {\em meshsize} the size of an individual finite element,
evaluated by means of an appropriate scalar measure (See below).
In the context of mesh adaptation,
the mesh generation algorithm is constrained
by a {\em meshsize field} defined on the domain to be meshed,
and whose value at a point
is the expected element size in the vicinity of that point.
The meshsize field is usually derived from an error estimation procedure
performed on the solution
of a prior finite element or finite volume analysis,
by requesting a smaller meshsizes
at places where the discretization error is deemed large.

Yet, when solving a problem for the first time,
an initial mesh has to be generated
without information from a prior computation,
and the meshsize field to generate that initial mesh
has to be constructed from crash
on basis of the geometrical data of the model only.
Given a CAD model, there exist a number of
theoretical prerequisites on the meshsize field 
to ensure a computable mesh,
and this is the purpose of this paper
to describe an automated algorithm
to compute a meshsize field fullfilling those prerequisites  {\it a priori}.
The proposed approach is ``user-driven'',
in the sense that users should be able to
generate a workable computational mesh in one click
on basis of a limited 
number of intuitive meshing parameters,
understandable by any finite element designer
with no extensive background in meshing.

An isotropic meshsize field  $h({\bf x})$
is thus a scalar function indicating the expected element size
at any point ${\bf x}$ in a domain to be meshed.
A first design choice 
concerns the mathematical representation of $h$.
As our goal is to build a ``first mesh'',
no background mesh is yet available
against which $h$ could be interpolated.
A classical solution (e.g. in Gmsh \cite{geuzaine2009gmsh})
is to define meshsizes 
directly on the geometrical entities of the model.
Meshsizes can be prescribed at the vertices of the CAD model, for instance,
and smoothly interpolated on model edges.
They are then subsequently interpolated 
on surface mesh vertices.
Meshsize fields interpolated this way
may however be biased by geometric features of the surface mesh,
such as gaps, fins or channels, which are assigned locally a small meshsize
that is not expected to spread out at distance in the bulk of the volume.
This approach is therefore not 100\% reliable
and defining the meshsize on auxiliary objects allows
for a better control and prevents the aforementioned phenomenon.

Two kinds of representation for meshsize fields
are encountered in the literature:
namely simplicial background meshes
\cite{cunha1997automatic, owen1997neighborhood,
  chen2017automaticSurf, chen2017automaticVol},
and Cartesian grids,
initially in the form of uniform grids \cite{pirzadeh1993structured},
and later on in the form of non-uniform 
octrees \cite{zhu2002background, tchon2005three, quadros2010computational}.
A graphical representation of such data-structures in the two-dimensional case
can be found in the first figure of \cite{persson2006mesh}.
It is immediately observed that uniform Cartesian grids are constrained by
the smallest feature in the CAD model. As that the number of
grid nodes grows cubically with the number of divisions,
the memory cost of uniform Cartesian grids quickly becomes prohibitive in practice,
and they were rapidly abandoned
for the sake of simplicial background meshes and tree-based grids.

Both 
representations have however their pros and cons.
Simplicial background meshes offer an accurate representation of
the boundaries, and the meshsize query procedure in the
three-dimensional case is reduced to a search
in the two-dimensional parametric space \cite{chen2017automaticSurf, chen2015automatic}.
However, meshsize fields represented this way
are rather sensitive to the location of the vertices in
the background mesh \cite{owen1997neighborhood},
and the access time to
the meshsize at a point in the background mesh
might be linear in the number of nodes in the worst case.
To improve on this,
Chen et al. \cite{chen2017automaticSurf} proposed a
\emph{walk-through} algorithm based on the backward search from Shan
et al. \cite{shan2008robust}, although it requires a well-guessed
element to quickly locate the point in the mesh, assuming an
interconnection between the background mesh and the meshing algorithm.

Octrees are orientation sensitive Cartesian structures.
They lack the geometrical flexibility of simplicial meshes
and significative
refinement may be necessary to accurately 
resolve surface-based
information, e.g. curvatures. On the other hand,
octree-based 
meshsize fields offer adaptive capabilities
to represent quickly and easily complex 
meshsize distributions across the structure.
Moreover, octrees offer 
fast access to any query points in $\mathcal{O}(\log_8 n)$ time,
where $n$ denotes the number of octants,
i.e., the number of leaves in the octree.
In this work, we use the 
octree implementation provided by \textsc{p4est}
\cite{BursteddeWilcoxGhattas11}. The serial version
of \textsc{p4est} is used, although
a scalable parallel implementation exists in the \textsc{p4est} library.
In our implementation, a uniform meshsize is assigned
to each octant in the octree,
which is the most natural option with \textsc{p4est}.

As mentioned earlier,
there are theoretical prerequisites
to ensure the computability of a mesh on a given CAD model.
Those prerequisites can be associated with
the following five intuitive mesh parameters.

\noindent {\bf Bulk size.} A bulk or default meshsize $h_b$. When
  creating the meshsize field, the octree is refined uniformly
  until every octant size is smaller or equal to $h_b$.
  This amounts to say that all meshsizes are initially set
  to the bulk value $h_b$.

\noindent{\bf Curvature.}
When using piecewise linear elements,
the main term of geometrical error produced by a mesh
is linked with the curvature of surfaces.
The local meshsizes $h({\bf x})$ should hence
be related to the maximal curvature
$\kappa_n({\bf x})$ of the surfaces.
This is done with the node density parameter $n_d$
that specifies the number of subdivisions
of the perimeter of the local osculating circle.

\noindent {\bf Small features.}
A CAD model may also contain narrow or thin regions, or \emph{features},
e.g. longerons that are the load-bearing components of
aerospace structures.
As such regions may have moderate or no curvature at all,
they are likely to be overlooked by an algorithm
that solely links meshsize with curvature.
The thickness of narrow regions can be estimated geometrically
thanks to the concept of medial axis, see Section 2.2.
On this basis a third parameter $n_g$ is defined in our algorithm
that specifies the minimal number 
of elements across the thickness of a narrow region.
We call \emph{feature meshsize } $h_f$ that thickness divided by $n_g$.
The feature meshsize $h_f$ and the curvature meshsize $h_c$
are the indicators used by our algorithm
to recursively refine the octree containing the meshsize field, see Section 2.3.

\noindent {\bf Boundedness.} A minimum mesh size $h_{\min}$ has to be
  defined as a fourth parameter in our algorithm
  to forbid inacceptably small meshsizes, whenever
  curvatures are very high for instance (e.g. at corners or at the tip of a cone).

\noindent {\bf Smoothness.} Accurate finite element and finite volume
  simulations usually require that meshsizes vary not too abruptly
  across the domain of computation.
  Yet, curvatures and  feature sizes
  may exhibit sharp variations in practice,
  resulting in unacceptably large meshsize gradations.
  A fifth parameter $\alpha > 1$ is thus defined that
  bounds (from above) the length ratio between two adjacent edges in the mesh.
In Section 2.4, we show that this
  condition is equivalent to limit the meshsize gradient to $|\nabla h| < \alpha-1$.

Our approach thus defines five parameters
that are \emph{easy to understand} by finite element practitioners,
and can be given a reasonable and rather \emph{universal} default value:

\begin{itemize}
	\item the bulk size, or default meshsize on newly created octants, $h_b$.
      The default value is $h_b = L/20$, where $L$ is the largest dimension
      of the axis-aligned bounding box of the CAD model;
	\item the minimal size allowed in the final mesh, $h_{\min}$.
      The default value is $h_{\min} = L/1000$;
	\item the number of elements $n_d$ used to accurately discretize
      a complete circle. The default value is $n_d = 20$;
	\item the number of element layers in thin gaps $n_g$.
      The default value is $n_g = 4$;
	\item the gradation, or length ratio of two adjacent edges
      in the final mesh, $\alpha$. The default value is $\alpha = 1.1$.
\end{itemize}

This meshsize field computation has been implemented in Gmsh,
which is also the tool used to generate all meshes presented in this paper.
It is planned that the algorithm presented in this paper be soon integrated
as a standard procedure in the meshing pipeline of Gmsh.


\section{Description of the algorithm : a worked-out example}
\begin{figure}[h!]
  \centering
	\includegraphics[width=0.24\textwidth]{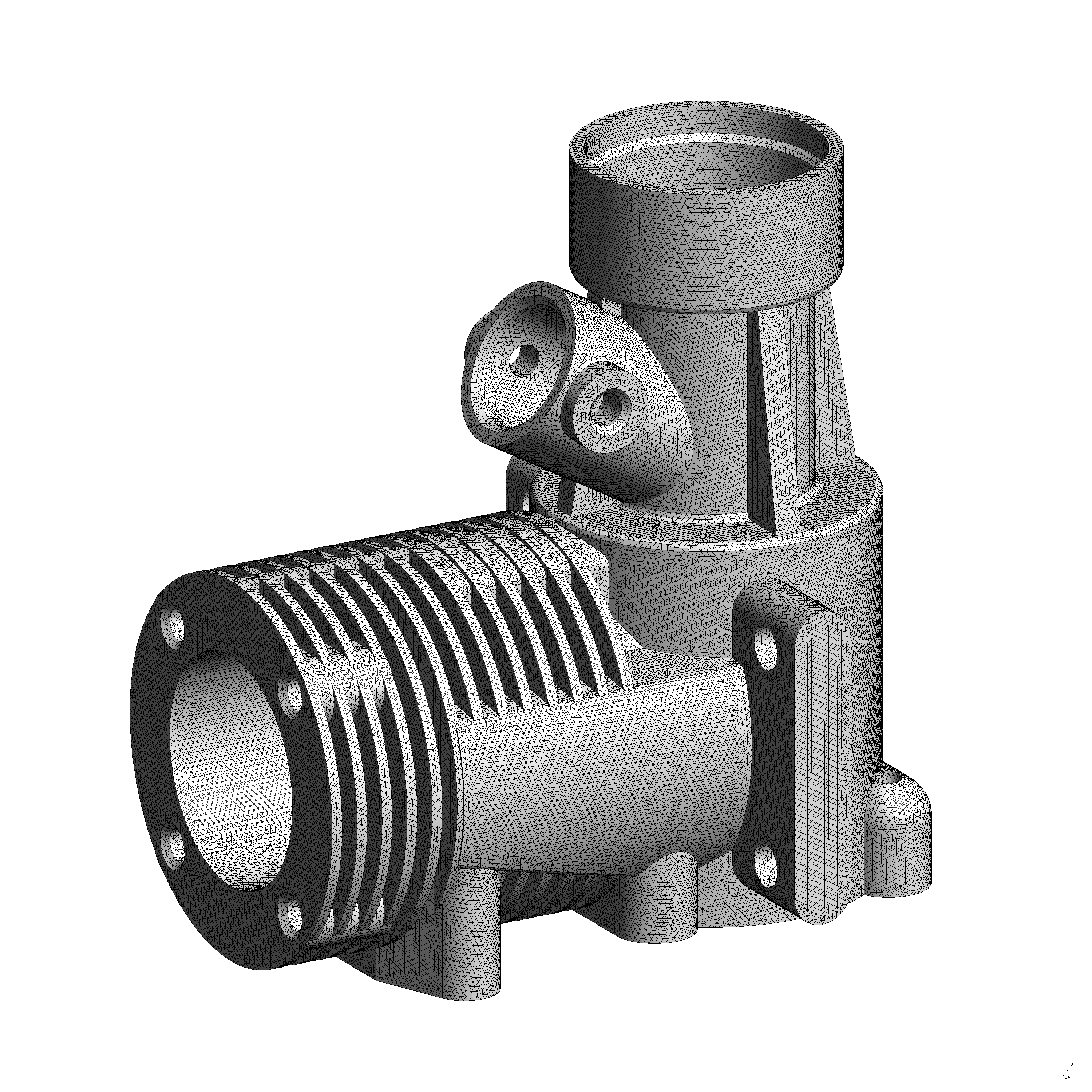}
	\includegraphics[width=0.24\textwidth]{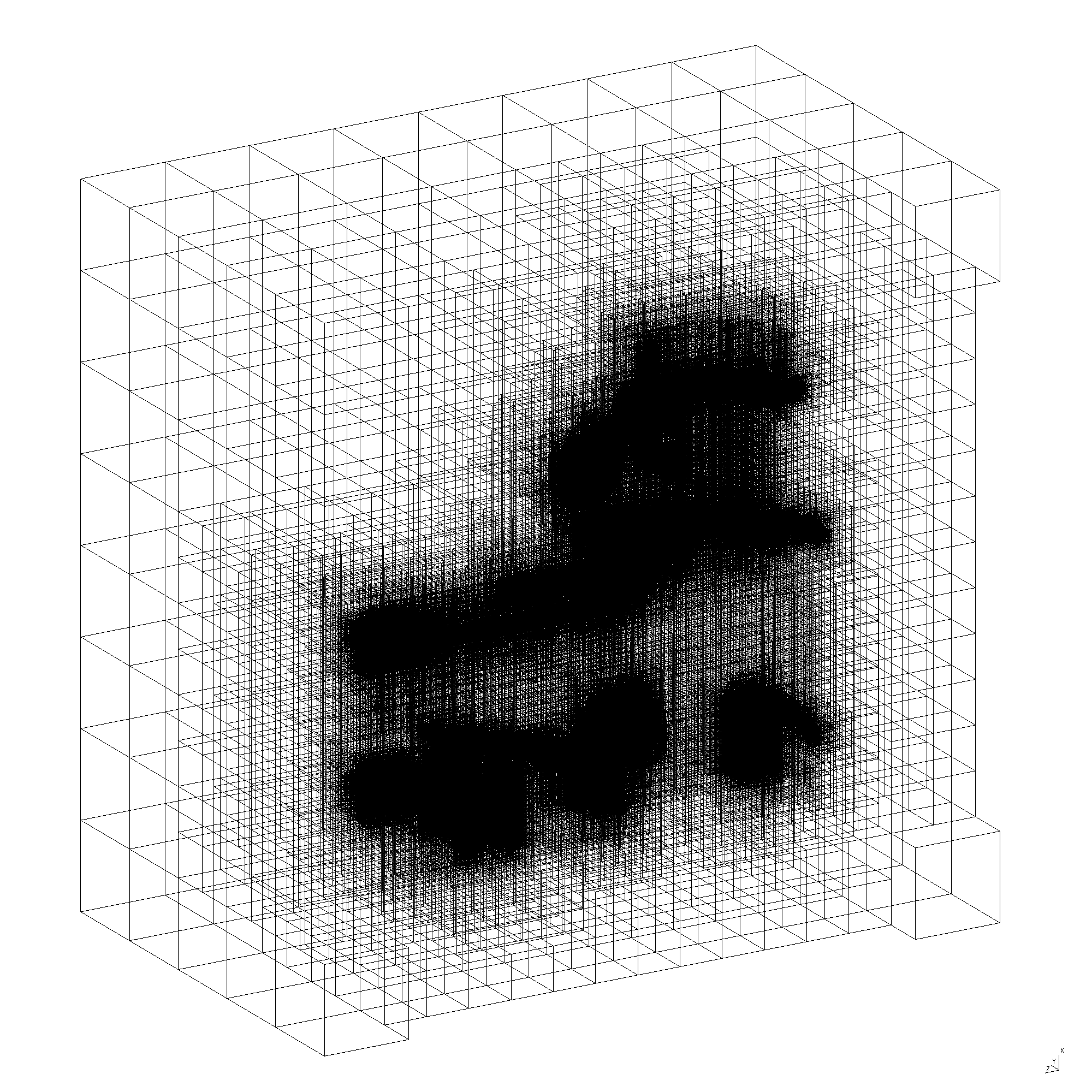}
	\includegraphics[width=0.24\textwidth]{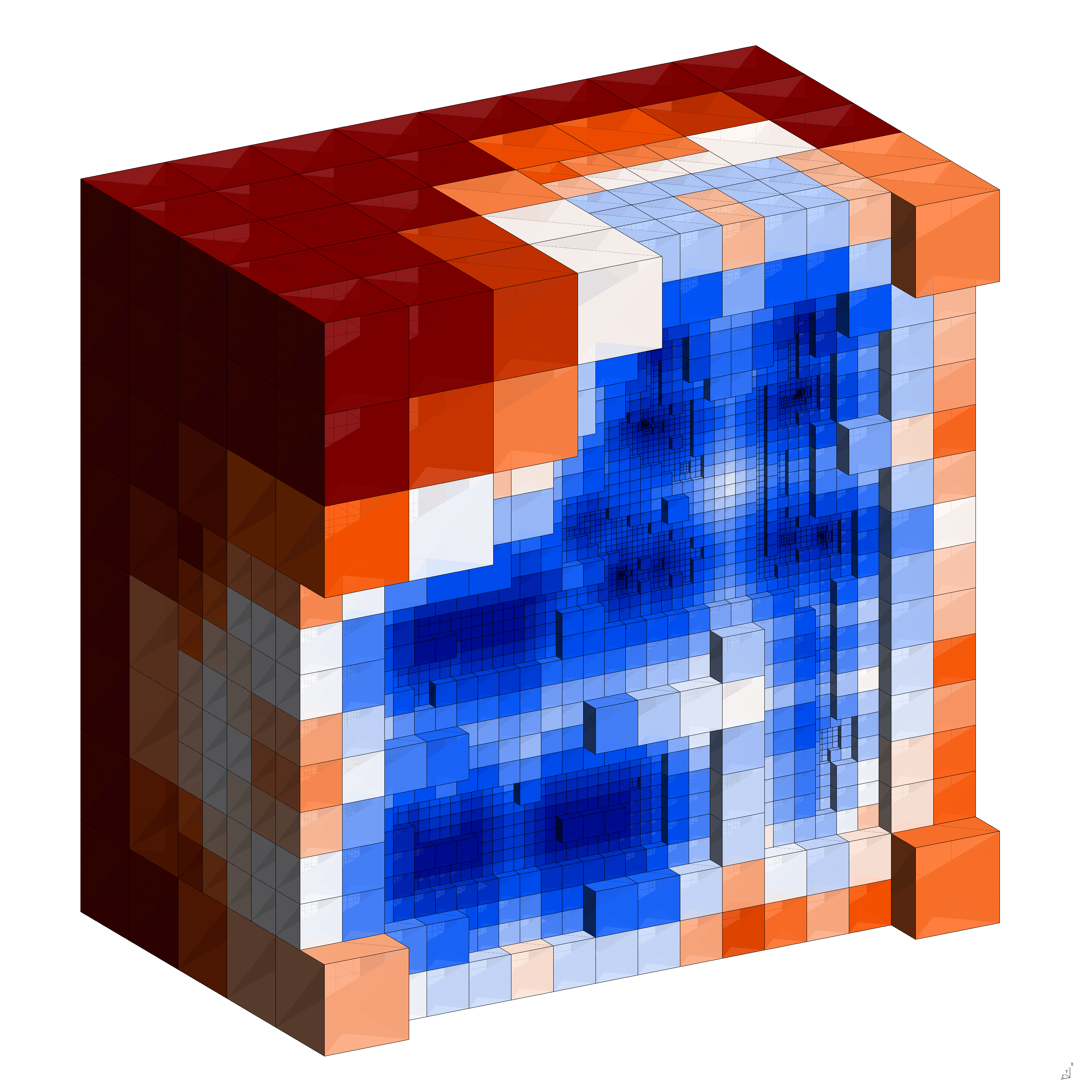}
	\includegraphics[width=0.24\textwidth]{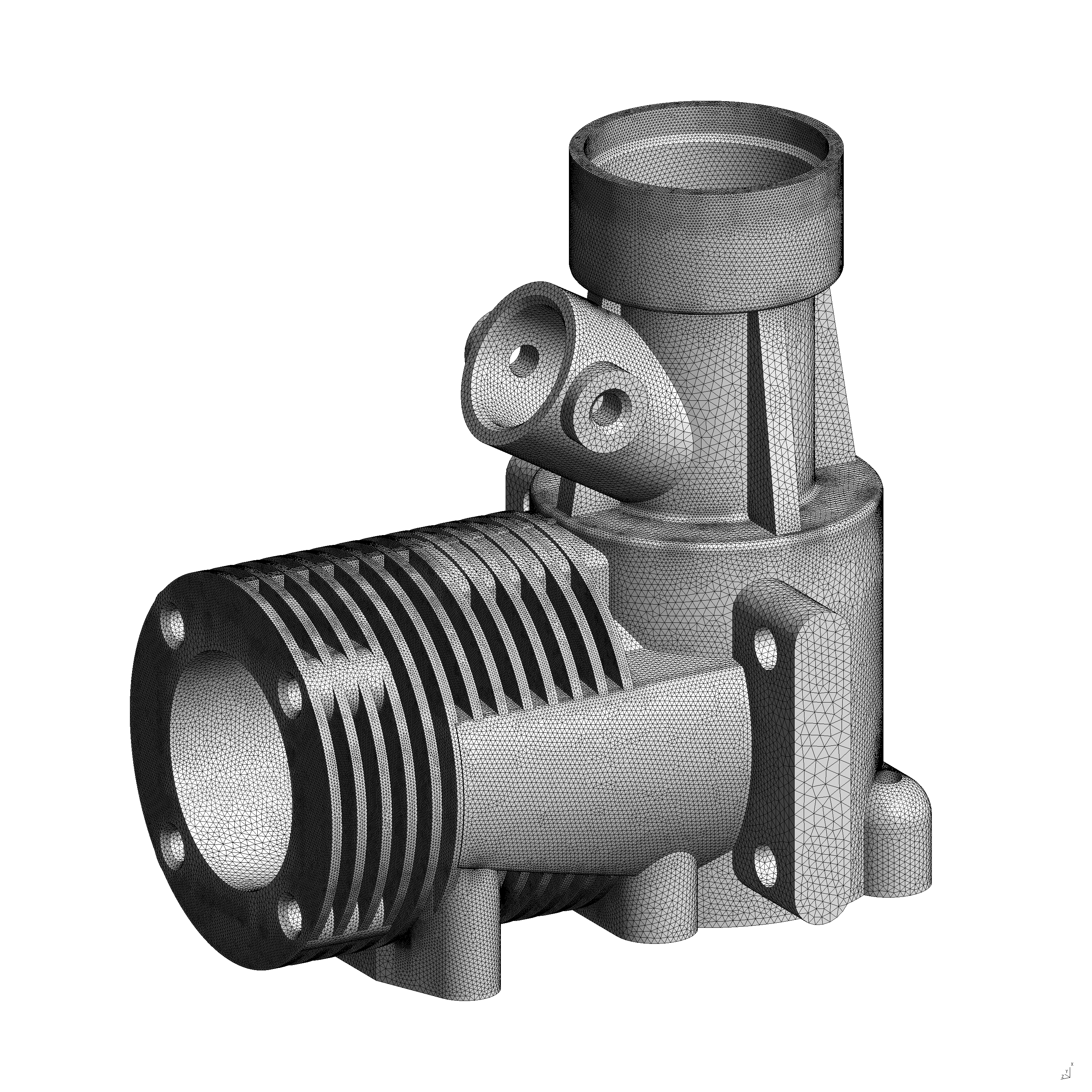}
	\caption{Overview of the algorithm for meshsize field computation, from left to right : $(i-ii)$ surface mesh of the engine block, from which discrete curvature and feature sizes are computed; $(iii-iv)$ the octree is refined based on curvature and feature meshsizes,
$(v)$ meshsize gradient is limited, yielding smoother meshsize field; generation of the final mesh.}
	\label{fig:geometry}
\end{figure}

In order to illustrate the steps of the construction of the
meshsize field $h({\bf x})$, the CAD model of an engine block is considered
as a application example
(Fig. \ref{fig:geometry}). This geometry
contains curved surfaces and narrow features that are
typical of real-life CAD models.
The input data for our algorithm is a surface mesh of the CAD model,
from which curvature and feature meshsize are computed.
The meshsize field is generated in an independent structure in the five following
steps (Fig. \ref{fig:geometry}) : $(i)$ compute the curvature meshsize $h_c$ from the approximate curvature on the surface mesh of the model;  $(ii)$ compute the feature meshsize $h_f$ from the medial axis of the geometry; $(iii)$ initialize the octree as the bounding box of the model and refine it uniformly until the size of
all octants is at most the bulk size $h_b$; $(iv)$ recursively refine the
octree based on both the curvature and the feature meshsize,
and assign the appropriate uniform meshsize in all newly created octants;
$(v)$ smooth out the meshsize field so as to limit its gradient to
$\alpha-1$.
During step $(v)$, the structure of the octree is not modified : only its stored meshsize field $h(\mathbf{x})$ is limited to satisfy $|\nabla h| < \alpha-1$, see Section 2.4.

\subsection{Approximation of surface curvatures}


With piecewise linear meshes,
the main term of geometrical error
is due to the curvature of surfaces.
Meshsize should thus be reduced in areas of high curvature.
We introduce to this end the \emph{curvature meshsize} $h_c(\mathbf{x})$.
Although definition slightly varies in the literature \cite{chen2017automaticSurf, chen2017automaticVol, turneraautomatic, deister2004fully}, they all rely on the  subdivision of the perimeter of local osculating circles.
The osculating circle at a point of a planar curve is the circle that best approximates the curve in the vicinity of the point, i.e., having same tangent and same curvature.
On a smooth surface, there is thus an osculating circle in every direction,
and the most critical curvature meshsize is related to the minimal radius of these circles, or reciprocally to the maximum normal curvature $\kappa_{n,\max}(\mathbf{x})$.
Hence the following definition for the curvature meshsize :
\begin{equation}
  h_c(\mathbf{x}) = \frac{2\pi r(\mathbf{x})}{n_d}
  = \frac{2\pi}{\kappa_{n,\max}(\mathbf{x}) \, n_d},
  \label{curvatureMeshsize}
\end{equation}
with $n_d$ a used-defined node density.
Whenever a CAD model is available in the background, surface parametrizations are at-hand and normal curvature can be obtained through the solid modeller's API. As far as our algorithm i concerned however, input data is a triangulation and normal curvatures are approximated following the tensor averaging methodology described by Rusinkiewicz \cite{rusinkiewicz2004estimating}, briefly recalled here with their notations.

Let $(\mathbf{u},\mathbf{v})$ denote an orthonormal basis in the tangent plane
at a point $\mathbf{x}$ of a smooth surface,
and $\boldsymbol{s}=(s_1,s_2)$ be an arbitrary direction in that plane.
The normal curvature at $\mathbf{x}$ in the direction $\boldsymbol{s}$ is given by
$$
\kappa_n(\mathbf{x}) = \II(\boldsymbol{s},\boldsymbol{s})
= \begin{pmatrix}s_1 & s_2 \end{pmatrix} \, \II \, \begin{pmatrix}s_1 \\ s_2 \end{pmatrix} = \begin{pmatrix}s_1 & s_2 \end{pmatrix} \begin{pmatrix}e & f\\ f & g\end{pmatrix} \begin{pmatrix}s_1 \\ s_2 \end{pmatrix}
$$
where $\II$ denotes the second fundamental form.
The eigenvalues $\kappa_1$ and $\kappa_2$ of the symmetric matrix $\II$,
known as the \emph{principal curvatures},
are the maximum and miminum values of normal curvature at $\mathbf{x}$.
Since we need $\kappa_{n,\max}(\mathbf{x}) = \max(\vert \kappa_1 \vert,\vert \kappa_2 \vert)$ to define the curvature meshsize \eqref{curvatureMeshsize},
our goal is thus to build an approximation of $\II$
at each vertex of the surface mesh.

The idea in \cite{rusinkiewicz2004estimating} is to compute $\II$ first on the triangles, and average them over adjacent triangles to obtain the needed per-vertex information. We start by computing per-vertex normal vectors $\mathbf{n}_i$ by averaging the normals of all faces adjacent to each vertex. On each triangle, an arbitrary orthonormal coordinate system $(\mathbf{u}_f,\mathbf{v}_f)$ is then defined. The components of the quadratic form  $\II$ in that basis read
$$
\II = \begin{pmatrix} \II(\bfu_f,\bfu_f) & \II(\bfu_f,\bfv_f) \\ \II(\bfu_f,\bfv_f) & \II(\bfv_f,\bfv_f) \end{pmatrix}.
$$
They can be evaluated as $\II(\bfu,\bfv) = L(\bfu)\cdot\bfv$
where the \emph{shape operator} $L(\boldsymbol{s}) = \nabla_{\boldsymbol{s}}\mathbf{n}$
is the directional derivative of the normal vector $\mathbf{n}$
along a direction $\boldsymbol{s}$ in the tangent plane.
One has now the following finite difference approximation
\begin{equation}
  \II(\boldsymbol{e}_0,\bfu_f) = L(\boldsymbol{e}_0)\cdot\bfu_f = \nabla_{\boldsymbol{e_0}}\mathbf{n} \cdot \bfu_f = (\mathbf{n}_2-\mathbf{n}_1) \cdot \bfu_f,
  \label{eq1}
\end{equation}
with $\boldsymbol{e}_0$ the edge from vertex $\mathbf{x}_1$ to vertex $\mathbf{x}_2$.

On the other hand, the edge $\boldsymbol{e}_0$,
regarded as a vector, admits a decomposition $\boldsymbol{e}_0 = (\boldsymbol{e}_0\cdot\bfu_f)\bfu_f + (\boldsymbol{e}_0\cdot\bfv_f)\bfv_f$ in that basis and, because the second fundamental form is a bilinear form, it follows that
\begin{equation}
  \II(\boldsymbol{e}_0,\bfu_f) = \II\big((\boldsymbol{e}_0\cdot\bfu_f)\bfu_f + (\boldsymbol{e}_0\cdot\bfv_f)\bfv_f, \bfu_f\big) = (\boldsymbol{e}_0\cdot\bfu_f)\,\II(\bfu_f,\bfu_f) + (\boldsymbol{e}_0\cdot\bfv_f)\,\II(\bfv_f,\bfu_f).
  \label{eq2}
\end{equation}
Combining \eqref{eq1} and \eqref{eq2}, and proceeding the same way for $\bfv_f$ yields the following pair of relationships between the components of $\II$ :
$$\II \begin{pmatrix} \boldsymbol{e}_0\cdot\bfu_f \\ \boldsymbol{e}_0\cdot\bfv_f \end{pmatrix} = \begin{pmatrix} (\mathbf{n}_2 - \mathbf{n}_1) \cdot \bfu_f \\ (\mathbf{n}_2 - \mathbf{n}_1) \cdot \bfv_f \end{pmatrix}.$$
Repeating the same procedure for the two remaining edges of the triangle, one ends up with a system of 6 equations for 3 unknowns, which can be solved using a least square method.

\begin{figure}[h!]
  \centering
  \begin{subfigure}{.4\textwidth}
  \includegraphics[width=\linewidth]{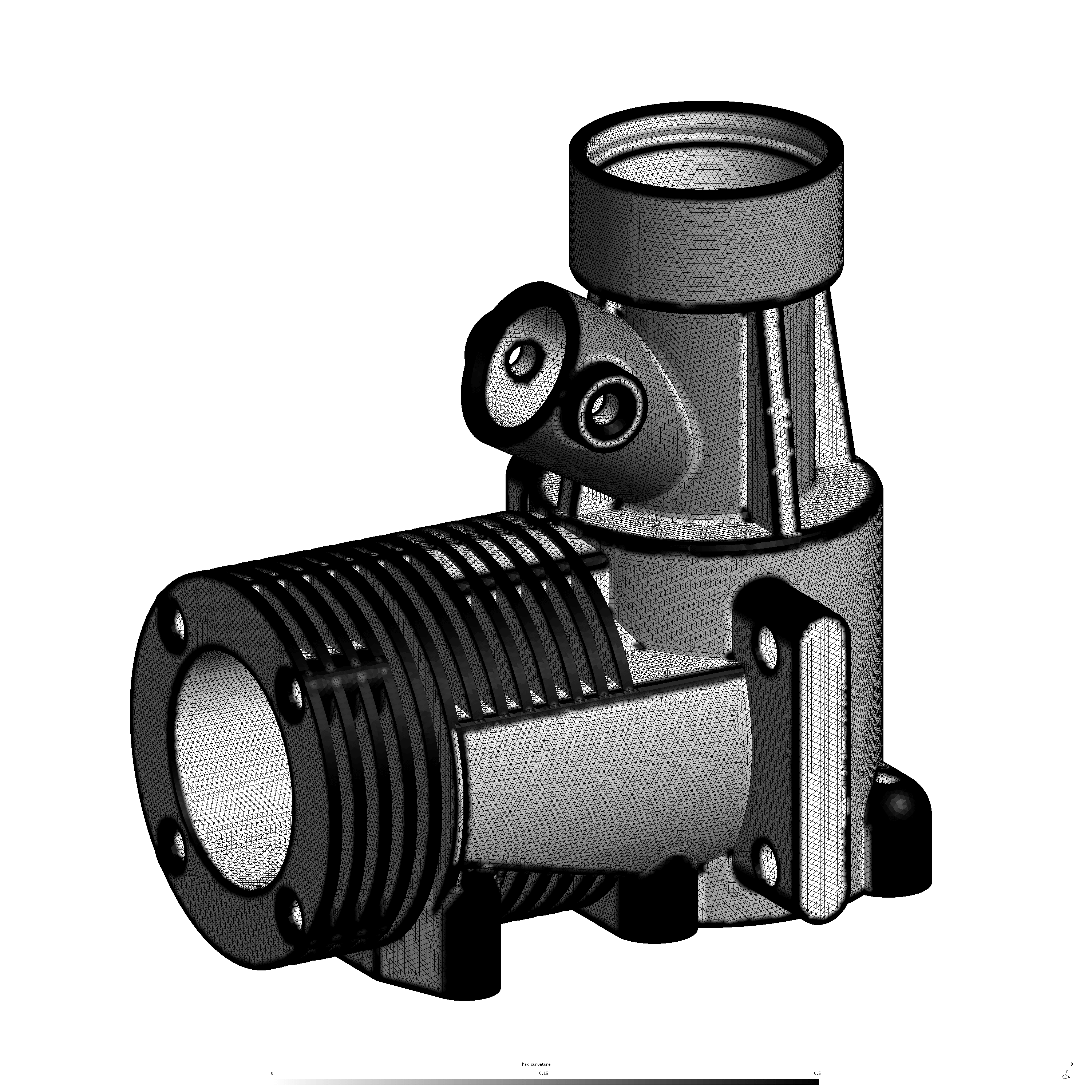}
  \caption{}
  \label{fig:curvature}
  \end{subfigure}
  \hspace{1.5cm}
  \begin{subfigure}{.4\textwidth}
  \includegraphics[width=\linewidth]{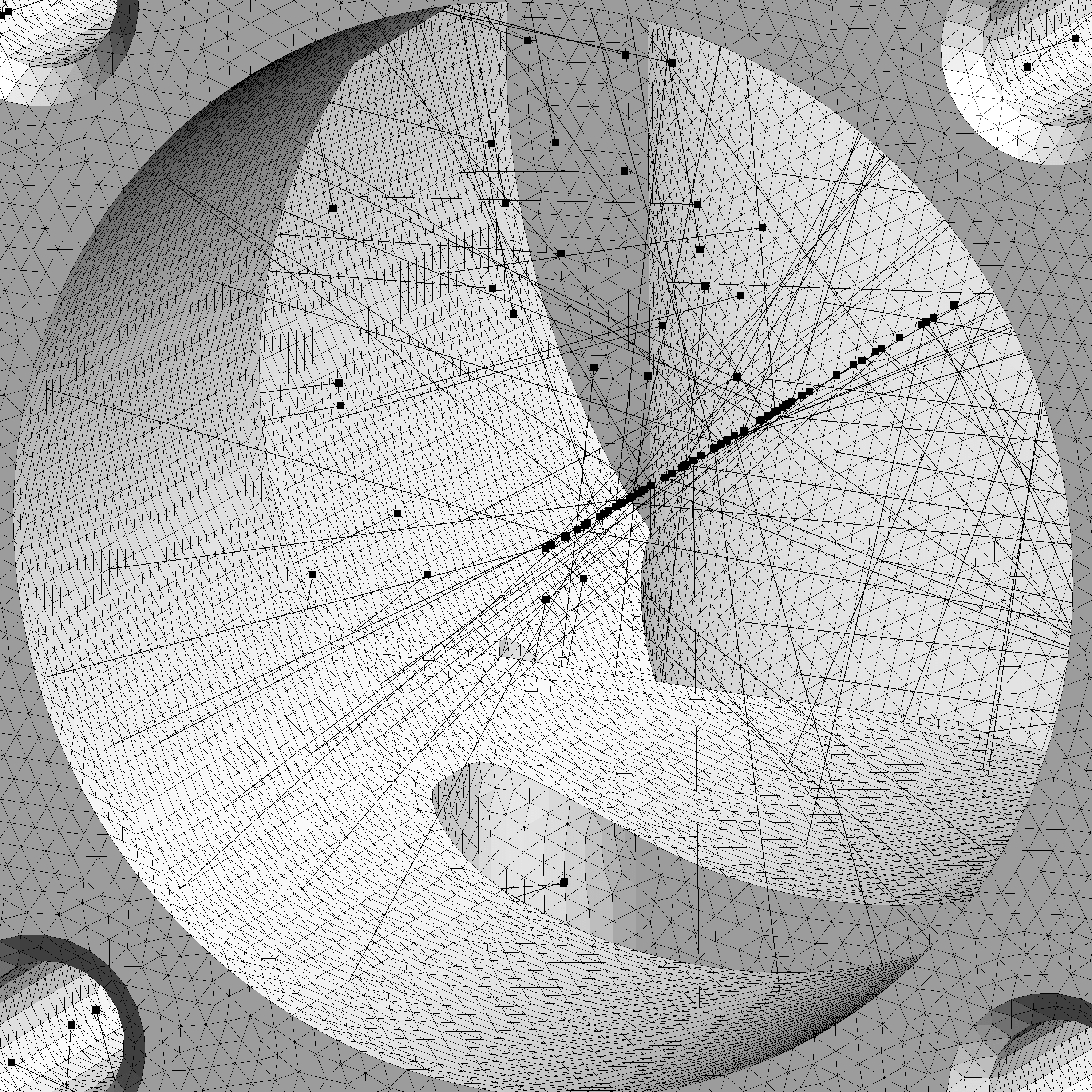}
  \caption{}
  \label{fig:poles}
  \end{subfigure}
  \caption{Left : Approximated maximum normal discrete curvature $\kappa_n = \max(\vert\kappa_1 \vert,\vert\kappa_2\vert)$ of the surface triangulation, shown in grayscale : curved areas are in black while regions with mild curvature are in light gray. Right : A subset of the Voronoï vertices lying the farthest from the Delaunay vertices, called \emph{poles} (black dots), approximate the medial axis as the mesh density increases.}
\end{figure}

This per-face approximation of $\II$ is expressed in the local basis $(\bfu_f,\bfv_f)$. In order to now combine the contributions of all triangles adjacent to a vertex $p$, one further orthonormal basis $(\bfu_p,\bfv_p)$ is defined in the plane perpendicular to the normal vector at $p$. The basis $(\bfu_f',\bfv_f')$ is defined as the $(\bfu_f,\bfv_f)$ basis slightly tilted (i.e., rotated) to be coplanar with $(\bfu_p,\bfv_p)$, and the contribution of the face to $\II_p$ can then be expressed in the rotated basis as
$$
e_p = \bfu_p^T\,\II\,\bfu_p = \begin{pmatrix} \bfu_p \cdot \bfu_f' \\ \bfu_p \cdot \bfv_f' \end{pmatrix}^T\II\,\begin{pmatrix} \bfu_p \cdot \bfu_f' \\ \bfu_p \cdot \bfv_f' \end{pmatrix}, ~~~~~~~ f_p = \bfu_p^T\,\II\,\bfv_p, ~~~~~~~ g_p = \bfv_p^T\,\II\,\bfv_p.
$$
As for the normals, the contributions of all triangles adjacent to the vertex $p$ are
then averaged to finally obtain  $\II_p$.
Finally, the curvature meshsize $h_c$ of the vertex is computed from the maximal eigenvalue of $\II_p$ using \eqref{curvatureMeshsize} (Fig. \ref{fig:curvature}).

\subsection{Feature size}
Whenever two surfaces with moderate or no curvature are close to each other but the distance between them is
smaller than the meshsize, the mesh generator will place only one element
in the gap between those surfaces. In many engineering applications like solid
mechanics or fluid mechanics, having only one element in a gap means
that both sides are connected by one single mesh edge. If
Dirichlet boundary conditions are applied,
such as a non-slip boundary condition, the gap is then
essentially closed, leading to an unwanted change of the
domain topology. Note that {\it a posteriori} error estimation will not detect large
errors in those closed gaps where the solution is essentially
constant.

\begin{figure}[h!]
\begin{subfigure}{.49\textwidth}
  \centering
  \begin{tikzpicture}
    \node[anchor=south west,inner sep=0] (image) at (0,0)
    {\includegraphics[width=1.2\linewidth]{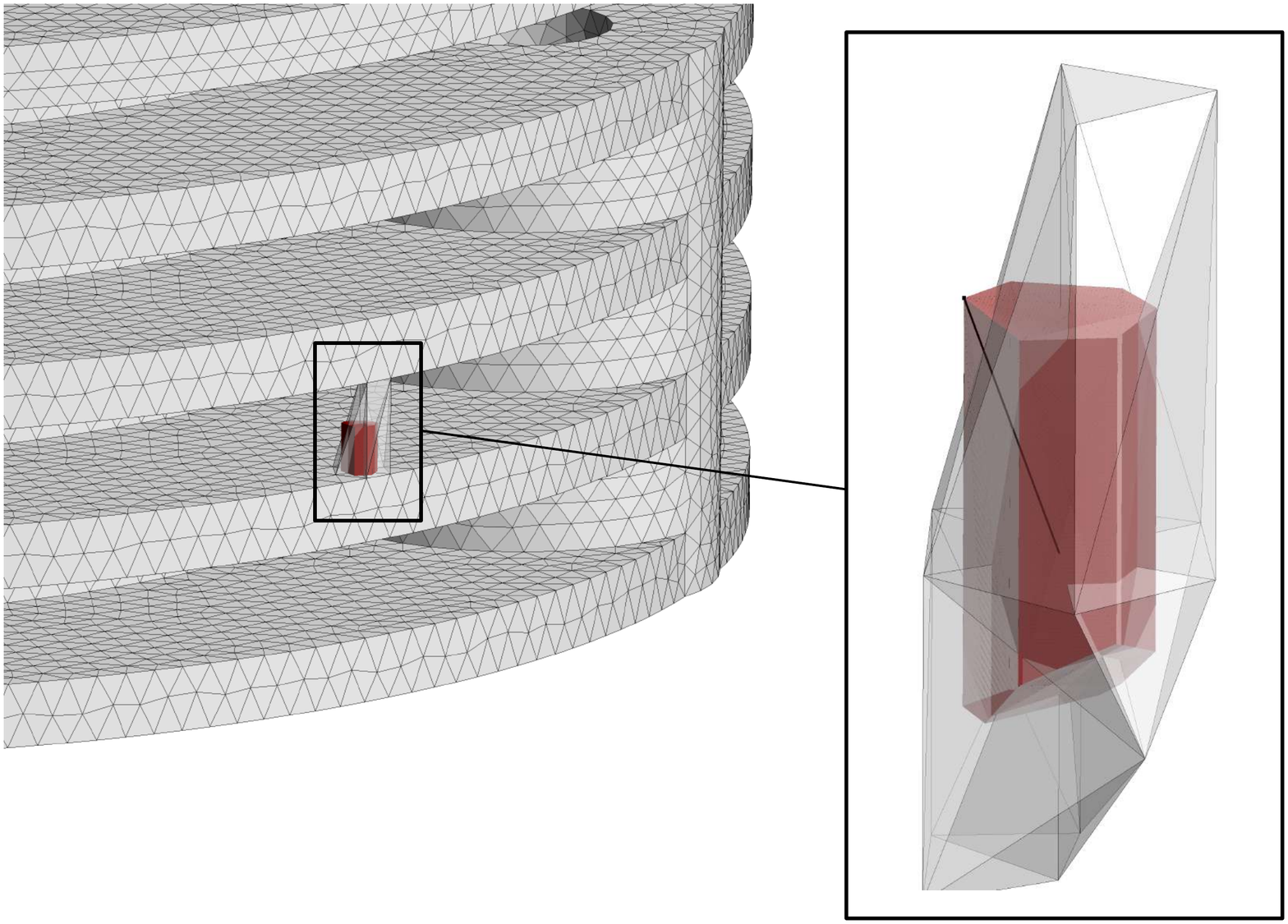}};
    \begin{scope}[x={(image.south east)},y={(image.north west)}]
      \node at (0.25,0.2) {$f(\mathbf{x})$};
      \draw[<->] (0.2,0.24) -- (0.2,0.31);
    \end{scope}
  \end{tikzpicture}
  \caption{}
  \label{fig:voronoiCell}
\end{subfigure}
\hspace{1cm}
\begin{subfigure}{.49\textwidth}
  \centering
  \begin{tikzpicture}
    \node[anchor=south west,inner sep=0] (image) at (0,0)
    {\includegraphics[width=0.6\textwidth]{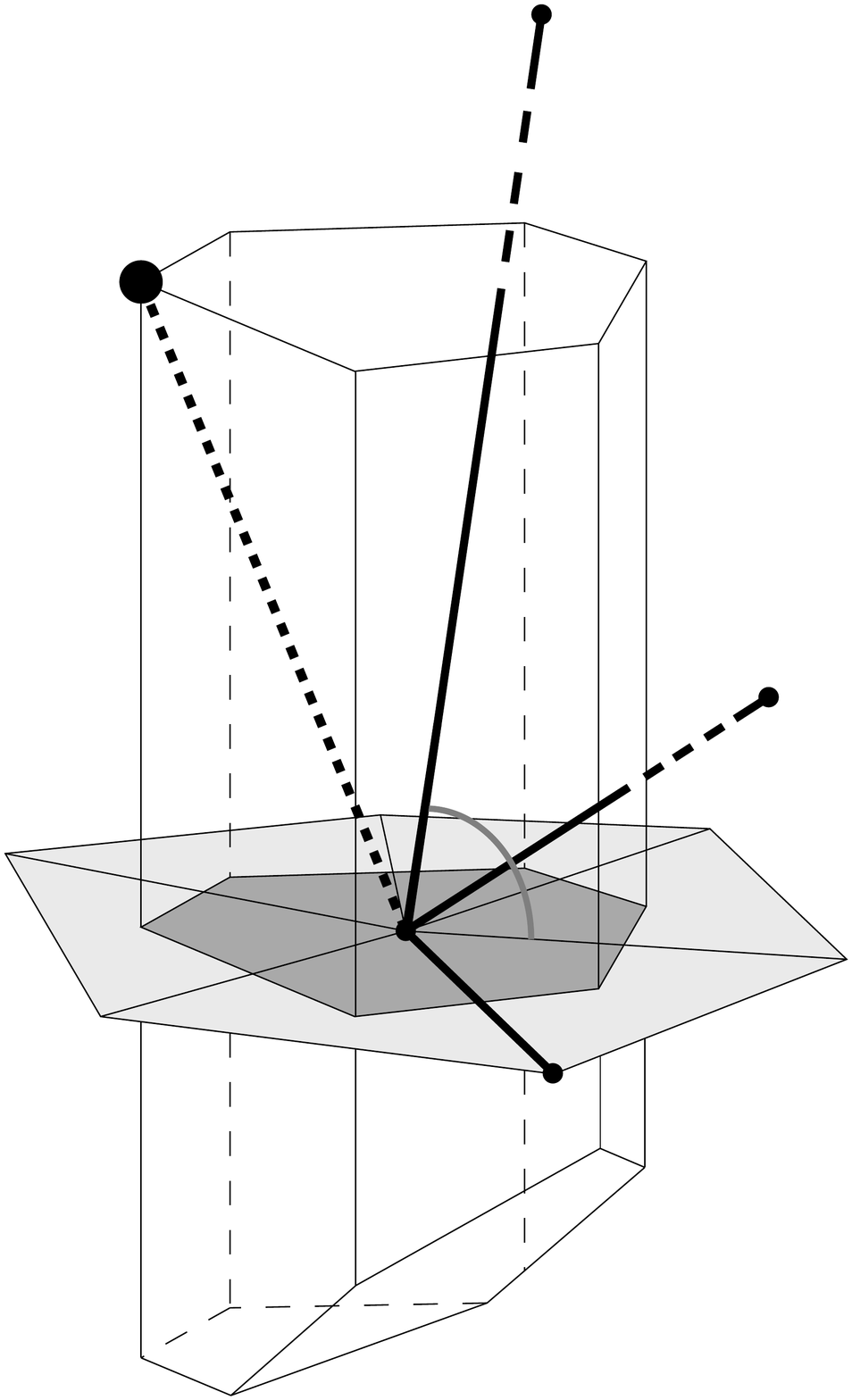}};
    \begin{scope}[x={(image.south east)},y={(image.north west)}]
      \node at (0.47,0.305) {$p$};
      \node at (0.56,0.99) {$q_1$};
      \node at (0.9,0.5) {$q_2$};
      \node at (0.65,0.27) {$q_3$};
      \node at (0.19,0.82) {$p^+$};
      \node at (0.37,0.6) {$v_p$};
      \node at (0.54,0.45) {$\theta$};
      \node at (0.95,0.35) {$U_p$};
    \end{scope}
  \end{tikzpicture}
  \caption{}
  \label{fig:voronoiCellSchema}
\end{subfigure}
\caption{Left : View of the input triangular surface mesh and of the Delaunay tetrahedra adjacent to a mesh vertex $p$. Zoom in on the same tetrahedra (transparent) and the associated dual Voronoï cell (in red). Right : Voronoï cell dual to the Delaunay vertex $p$ in the triangular surface mesh. The pole vector $v_p$ (dashed)
connects the mesh vertex $p$ to the farthest Voronoï vertex (tetrahedron circumcenter) or \emph{pole} $p^+$. The Delaunay edges $pq_1, pq_2$ and $pq_3$ connect $p$ to neighbouring Delaunay vertices $q_1, q_2$ and $q_3$, respectively. Only $pq_1$ and $pq_2$ satisfy one of the filtering conditions and are considered
to compute the feature size.
In grey, the umbrella $U_p$ of $p$ is the Delaunay facets (triangles) dual to the Voronoï edges cut by the plane through $p$ with normal $v_p$ (not shown).}
\end{figure}

Hence, special care should be given to such narrow geometrical
regions, where curvature information alone is not enough to determine a suitable meshsize field.
To this end, we define the \emph{feature size} $f(\mathbf{x})$ as a measure of the local gap
thickness (Fig. \ref{fig:voronoiCell}).
If $\partial V$ is the boundary of a volume $V$,
$f(\mathbf{x})$, $\mathbf{x} \in \partial V$,
is defined as twice the distance between $\mathbf{x}$ and the \emph{medial axis} of volume $V$.
If a surface 
bounds two volumes, then the
minimum feature size is chosen.
The \emph{feature meshsize} at the considered vertice, now, is the feature size devided by the desired number of element layers in narrow regions $n_g$, i.e.,  $h_f(\mathbf{x}) = f(\mathbf{x})/n_g$. The feature meshsize $h_f$ is thus in a similar relationship to the feature size $f$ than the curvature meshsize $h_c$ was to the maximum principal curvature $\kappa$.

The evaluation of the feature size $f(\mathbf{x})$ thus requires computing an approximation of the medial axis of all volumes in the computational domain.
The medial axis of a volume $V$, also referred to as its \emph{skeleton}, is defined as the set of points having more than one closest point on its boundary $\partial V$. Equivalently, the medial axis is the set of the centers of all spheres tangent to $\partial V$ in two or more points, and the feature size $f(\mathbf{x})$ is twice the radius of the sphere tangent at $\mathbf{x}$. 

We rely on the algorithm \textsc{Medial} introduced by Dey and Zhao \cite{dey2004approximating}
to compute a discrete approximation of the medial axis.
The algorithm is based on the
\emph{Voronoï diagram of the vertices of the surface mesh},
and it has suitable convergence properties in the sense that the output set of facets converges to the medial axis as the surface mesh density increases.
It takes as input data a Delaunay tetrahedrization of the vertices \emph{of the surface mesh only},
i.e., a set of tetrahedra, often called \emph{empty mesh}, filling the computational domain and whose nodes all lie on the surfaces.
The rationale behind the algorithm \textsc{Medial} can be sketched as follows
with the notations of Dey and Zhao.
In 2D, the vertices of the Voronoï cells that are  dual of an empty mesh
give straigth away an approximation of the medial axis. The same is however not always true in 3D,
because sliver tetrahedra of the empty mesh 
can persist close to the boundary as the surface mesh is refined. However, by pruning as explained below the Voronoï vertices dual of these sliver tetrahedra, a subset of Voronoï vertices called \emph{poles} can be defined
that do approximate the medial axis (Fig.~\ref{fig:poles}).
Given the Voronoï cell 
dual of a surface vertex $p$ (Fig.~\ref{fig:voronoiCell}), the corresponding pole $p^+$ is defined as the Voronoï vertex that is the most distant from $p$. Each 
surface vertex $p$ is thus associated with a pole $p^+$, and then with a \emph{pole vector} $v_p = p^+-p$ that approximates the normal to 
$\partial V$ at $p$ (Fig.~\ref{fig:voronoiCellSchema}).

The plane passing through $p$ with normal $v_p$ intersects edges of the Voronoï cell, and the Delaunay facets 
dual to these edges 
are pictorially called the \emph{umbrella} $U_p$ \emph{of} $p$
by Dey and Zhao (Fig.~\ref{fig:voronoiCellSchema}, in light grey).

Triangles of the umbrella are used to select some of the Delaunay edges
adjacent to $p$ whose dual facets will eventually form the discrete medial axis.
The idea is to select Delaunay edges $pq$ that $(a)$ make a sufficiently large angle with the triangles of the umbrella, or $(b)$ are significantly longer than the circumradius of these triangles. Condition $(a)$ measures how normal the edge $pq$ is to the umbrella $U_p$ and is referred to as the \emph{angle condition}. The edge should make an angle larger than $\theta$ with each triangle of the umbrella, or conversely, an angle smaller than $\pi/2-\theta$ with their normal vector. In practice, we evaluate \cite{dey2004approximating} : $$\max_{i\,\in\,U_p} \angle(\mathbf{pq}, \hat{\mathbf{n}}_i) < \frac{\pi}{2} - \theta,$$ where $\mathbf{pq}$ is a vector parallel to the edge $pq$, $i$ denotes a triangle in $U_p$ and the threshold angle is set to $\theta = \pi/8$. Condition $(b)$ ensures that edges of the surface mesh are removed. It does so by selecting long edges for which the angle condition has failed, and is referred to as the \emph{ratio condition}. In practice, only edges at least $\rho = 8$ times longer than the circumradius $R$ of the triangles in $U_p$ are considered, and the condition reads \cite{dey2004approximating} : $$\max_{i\,\in\,U_p}\frac{\Vert pq \Vert}{R_i} > \rho.$$
The numerical values for parameters $\theta$ and $\rho$ are those suggested by Dey and Zhao. They yield a good approximation of the medial axis for a large variety of geometries.
If a Delaunay edge $pq$ satisfies either of these two conditions, its dual Voronoï facet is added to a set $F$, forming the approximate medial axis. Consider the three mesh vertices $q_i$ connected to $p$ through a Delaunay edge $pq_i = q_i-p$ on
Figure~\ref{fig:voronoiCellSchema}. Edge $pq_1$ makes a large angle $\theta$ with the triangles of the umbrella, and edge $pq_2$ is several times longer 
 than the largest circumradius : both are added to the set $E$ dual to $F$.
Edge $pq_3$ is a short surface edge mesh lying flat to the umbrella, and is thus removed from the list of candidate edges.

In our meshsize field computation, we need not compute the dual facet to edges in $E$ : for each Delaunay edge $pq$ satisfying the angle or ratio conditions, the local feature size is directly given by the edge length $\Vert pq \Vert$. Meshsize at both vertices $p$ and $q$ is thus defined as the edge length divided by the
desired number of element layers in features :
$$h_f(\mathbf{x}_p) = h_f(\mathbf{x}_q) = \frac{\Vert pq \Vert}{n_g}.$$
The meshsize on the octants containing these vertices is then lower bounded if necessary :
\begin{equation*}
	h = \max\left( h_{\min}, \min \left( h_f, h_c, h_u, h_{b}\right) \right).
\end{equation*}

\begin{figure}[h!]
  \centering
  \begin{subfigure}{.4\textwidth}
    \centering
    \begin{tikzpicture}
      \draw (0,0) rectangle (7,3);
      \draw[very thick] (1.5,1.5) -- (5.5,1.5);

      \draw[dashed] (1.5,1.5) -- (0,3);
      \draw[dashed] (1.5,1.5) -- (0,0);
      \draw[dashed] (5.5,1.5) -- (7,3);
      \draw[dashed] (5.5,1.5) -- (7,0);

      \draw (1.5,1.5) circle (1.5cm);
      \draw (3.5,1.5) circle (1.5cm);
      \draw (5.5,1.5) circle (1.5cm);

      \draw (0.5,0.5) circle (0.5cm);
      \draw (0.25,0.25) circle (0.25cm);
      \draw (0.125,0.125) circle (0.125cm);
      \draw (0.0625,0.0625) circle (0.0625cm);

      \draw (0.5,2.5) circle (0.5cm);
      \draw (0.25,2.75) circle (0.25cm);
      \draw (0.125,2.875) circle (0.125cm);
      \draw (0.0625,2.9375) circle (0.0625cm);

      \draw (6.5,2.5) circle (0.5cm);
      \draw (6.75,2.75) circle (0.25cm);
      \draw (6.875,2.875) circle (0.125cm);
      \draw (6.9375,2.9375) circle (0.0625cm);

      \draw (6.5,0.5) circle (0.5cm);
      \draw (6.75,0.25) circle (0.25cm);
      \draw (6.875,0.125) circle (0.125cm);
      \draw (6.9375,0.0625) circle (0.0625cm);

      \draw[{Latex[length=2mm, width=2mm]}-{Latex[length=2mm, width=2mm]}] (3.5,0) -- (3.5,3);

      \draw[-{Latex[length=1mm, width=1mm]}] (5.5,3.7) -- (1.35,1.65);
      \draw[-{Latex[length=1mm, width=1mm]}] (5.5,3.7) -- (5.65,1.65);
      \node[above, text width=3.5cm,align=center] at (5.5,3.7) {Irrelevant branches in angles and corners};
      \node[right] at (3.5,2.5) {$f$};
    \end{tikzpicture}
  \caption{}
  \label{fig:secondaryBranches}
  \end{subfigure}
  \hspace{1.5cm}
  \begin{subfigure}{.49\textwidth}
  \includegraphics[width=\linewidth]{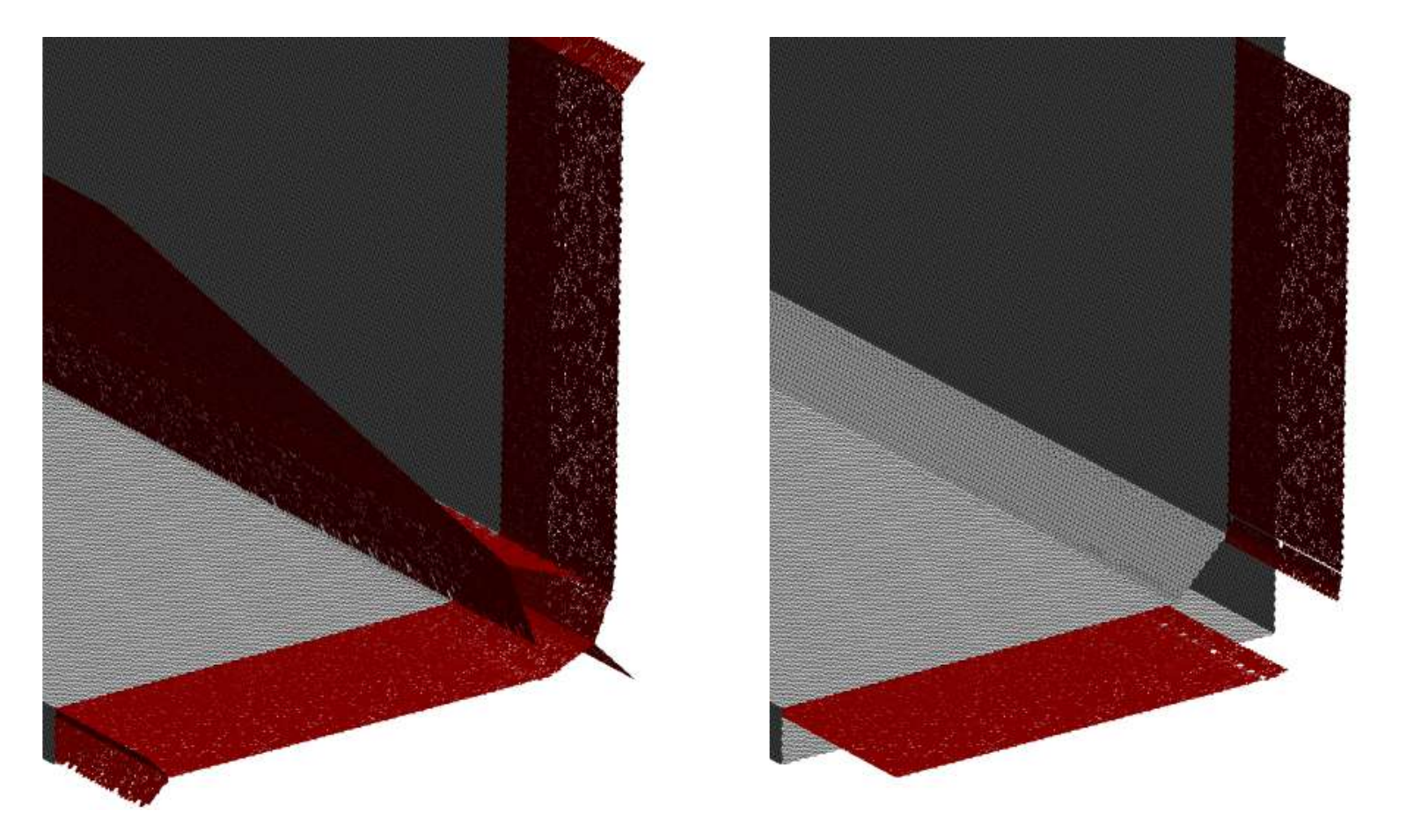}
  \caption{}
  \label{fig:wedge}
  \end{subfigure}
  \caption{Left : The medial axis of a rectangle consists of a main branch (thick) and four secondary branches (dashed) connecting the main branch to the corners. The radii of the spheres whose center lies on the main branch are representative of the feature size $f$ (here, the thickness of the rectangle), while the radii of the spheres on the secondary branches shrink to zero as the branch approaches the surface, and are thus not representative of $f$. To remove those secondary branches (Voronoï facets), the dual Delaunay edges are applied a second filtering process ensuring that the angle between the edge and the normal vectors at its extremities does not exceed $\theta$. Right : the medial axis (in red) of an angle geometry before filtering the branches in the corners (left) and the medial axis after filtering (right).}
\end{figure}
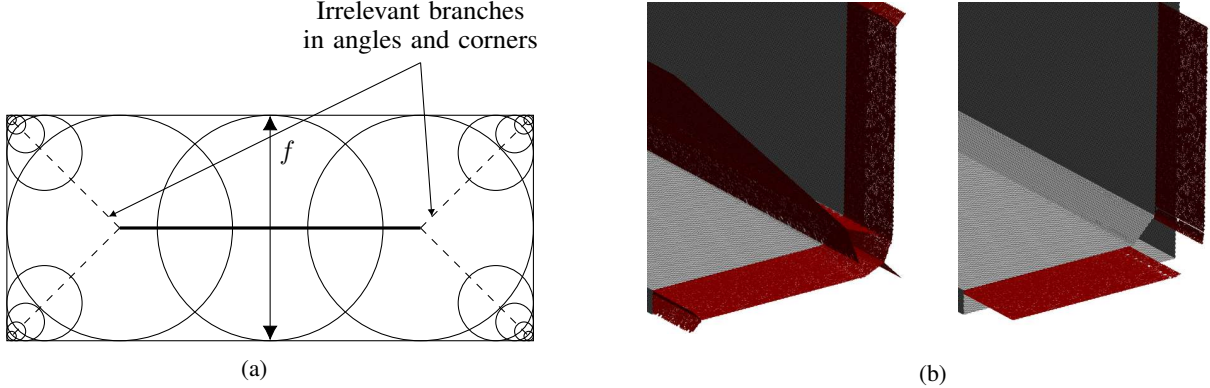

The medial axis also contains the centers of spheres with radius vanishing to zero in all angles and corners of the volume
(Fig.~\ref{fig:secondaryBranches}). Dual edges in these corners are smaller than the feature size one wants to identify, and they should thus be disregarded to avoid spurious small meshsize in such areas.
In practice, an edge $pq$ making an angle larger than $\theta$ with either one of the normal vectors $\hat{\mathbf{n}}_p$ and $\hat{\mathbf{n}}_q$ at its ends is also filtered out (Fig. \ref{fig:wedge}).

The quality of the approximation of the medial axis is directly related to the node density of the surface mesh. As pointed out by Dey and Zhao, the surface mesh should be an $\varepsilon$-sample, i.e., vertices should have neighbours within a distance $\varepsilon f(\mathbf{x})$, where $f(\mathbf{x})$ is the feature size and $\varepsilon$ is small. Of course, as the aim is here precisely to compute the feature size, it is not known beforehand whether or not the input surface mesh is a $\varepsilon$-sample. A first solution is to measure 
beforehand the most critical feature size $f_{crit} = \min_{\mathbf{x}} f(\mathbf{x})$ of the CAD model, then generate a uniform mesh with constant meshsize $h_{crit} = \varepsilon f_{crit}$, typically with $\varepsilon \leq 0.25$ as suggested in \cite{dey2004approximating}. The characteristic size of the input surface mesh is then constrained by the smallest feature in the geometry, which may result in an expensive size field computation.
While we could certainly adjust the meshsize of the input mesh to be $h_{crit}$ only in the small features, this would amount to manually specify the meshsize field, which we want to avoid. To circumvent this, one can compute an initial meshsize field based on a reasonably fine uniform mesh, then perform the mesh generation from this field. This intermediary mesh will include an initial refinement in the small features, and will be a better candidate for the final meshsize field computation.

\begin{figure}[h!]
  \centering
  \begin{subfigure}{.49\textwidth}
  \includegraphics[width=\linewidth]{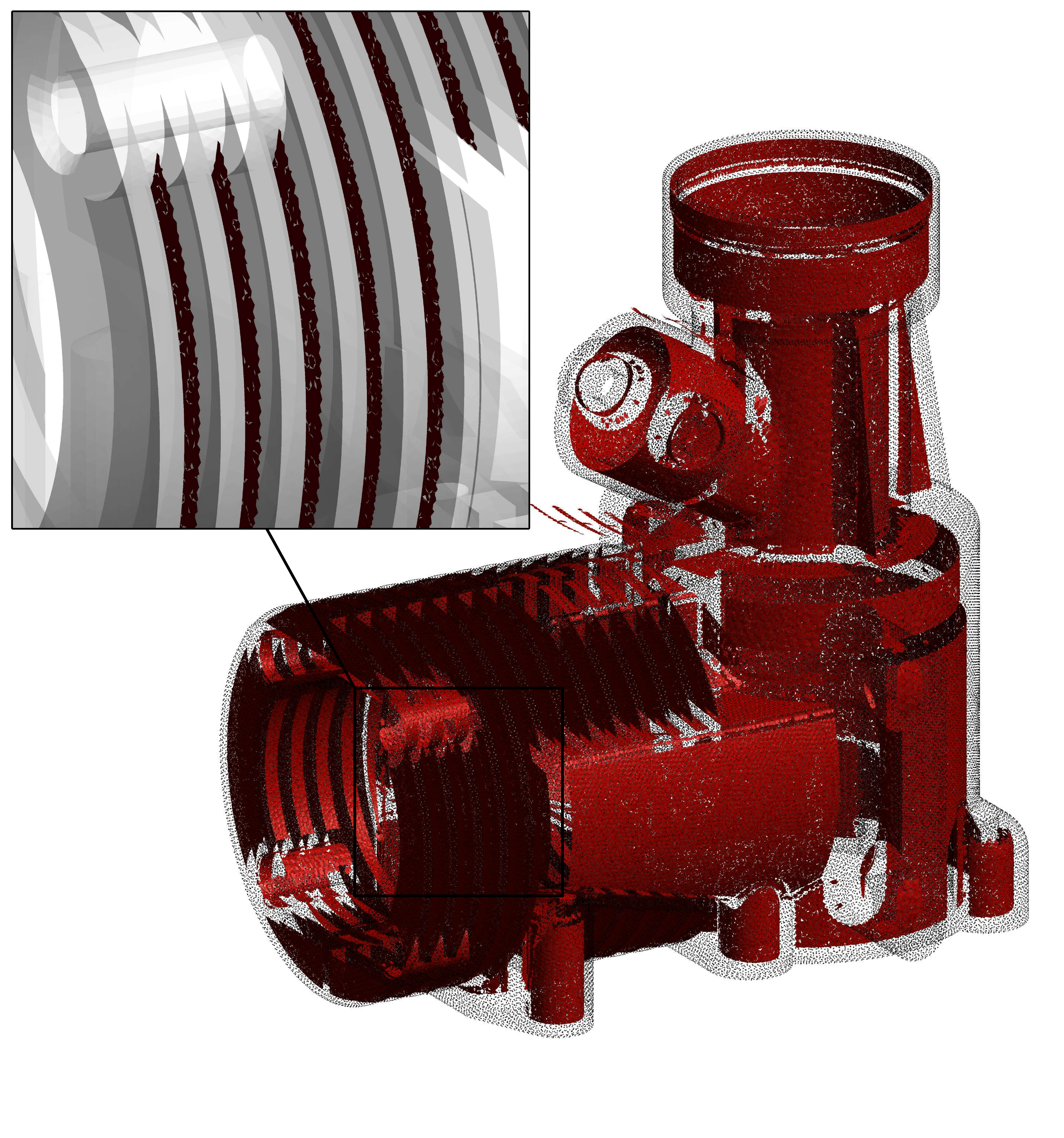}
  \caption{}
  \label{fig:medialAxis}
  \end{subfigure}
  \begin{subfigure}{.49\textwidth}
  \includegraphics[width=\linewidth]{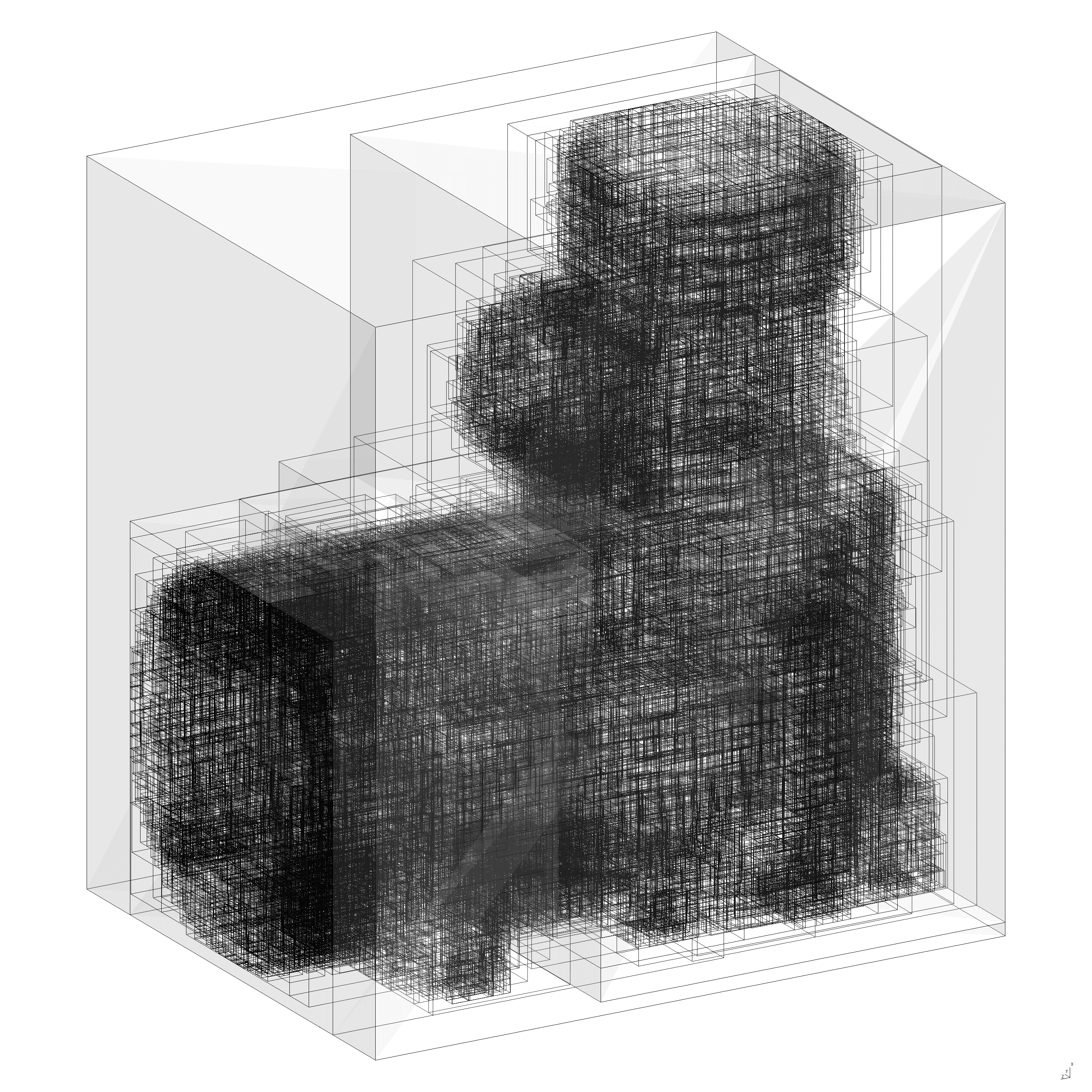}
  \caption{}
  \label{fig:rtree}
  \end{subfigure}
  \caption{Left : Approximated medial axis of the CAD model. The Voronoï facets dual to the filtered Delaunay edges are drawn in red. In the close-up view, only the facets of the medial axis lying outside of the volume are shown. Right: the R-Tree structure built from the bounding boxes of the triangles of the surface mesh.}
\end{figure}

\begin{figure}[h!]
	\includegraphics[width=0.5\textwidth]{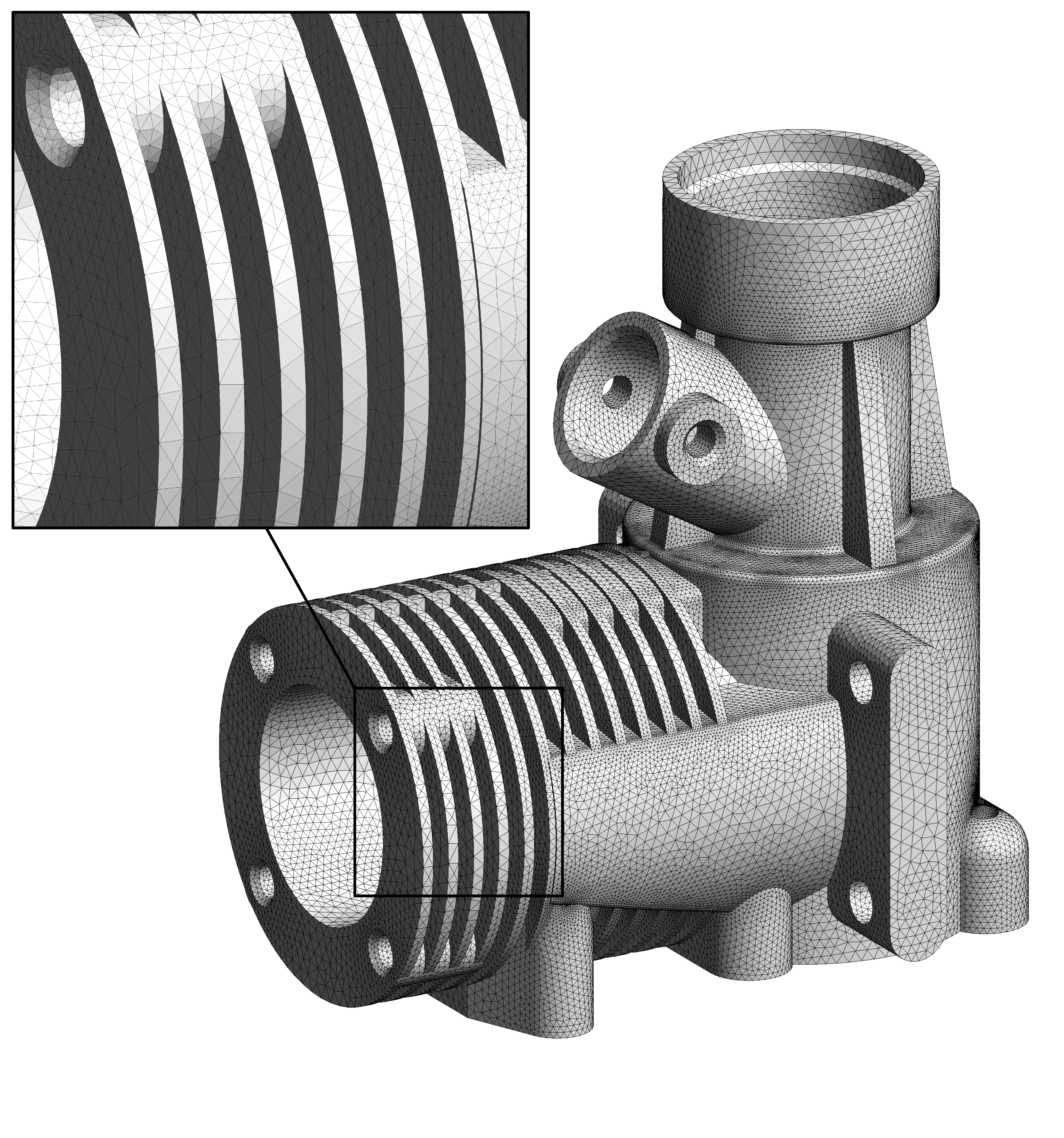}
  \includegraphics[width=0.5\textwidth]{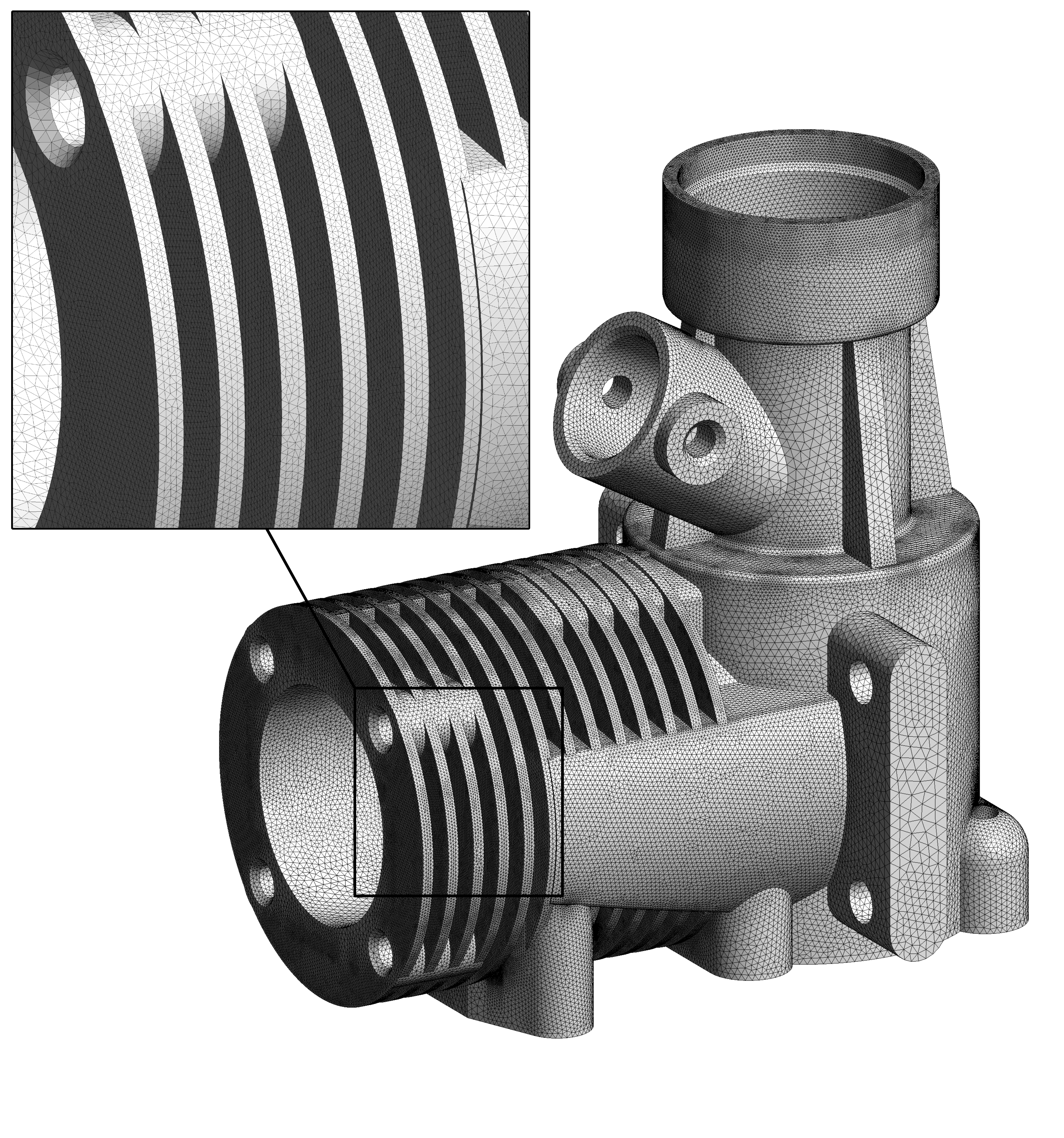}
  \caption{Surface mesh generated from the computed meshsize field : from curvature only (left) and considering both curvature and feature sizes (right). On the right, four layers of elements are generated in the fins around the largest cylinder.}
	\label{fig:withMedialAxis}
\end{figure}

The resulting medial axis for our block example after the two filtering operations is shown on Figure \ref{fig:medialAxis}. Taking into account the local feature size allows for a refined mesh in small features of the geometry, especially in areas with zero curvature, which would be overlooked otherwise (Fig. \ref{fig:withMedialAxis}, see also Fig. \ref{fig:blockVolumeMesh} for the volume mesh).

\subsection{Octree initialization and refinement}
The meshsize field $h(\mathbf{x})$ over the domain to be meshed
is now built as an octree structure.
The initial octant is defined as the
axis-aligned bounding box of the surface mesh,
stretched in all three dimensions by a factor $1.5$.
In order to ensure a suitable gradation in the mesh,
the main idea in leading the refinement of the octree is
that the dimension of each octant should eventually be
representative of the local meshsize.
The octree is then first subdivided recursively
and uniformly until the size of each octant is at most the bulk size $h_{b}$,
and the local meshsize value affected to the octants created during this initial step are set to $h_b$.
The octree is then further refined on basis of the curvature meshsize $h_c$ and the feature meshsize $h_f$. This information, which is available on the surfaces, has to be transfered to the octree, which is three-dimensional structure. One needs for that to detect efficiently the intersections between octants and the surface mesh.
To this end, the three-dimensional bounding box of each triangle of the surface mesh is added to an R-Tree \cite{beckmann1990r}, a data structure used for spatial access methods (Fig. \ref{fig:rtree}) which acts here as the intermediary between the geometry and the octree. More specifically, the R-Tree provides for each octant of index $i$ a list $\mathcal{T}_{i}$ of triangles whose bounding box intersects the octant. The octant with index $i$ is then divided until it becomes smaller than the minimal meshsize ($h_c$ or $h_f$) at the vertices of all triangles in $\mathcal{T}_i$. Whenever an additional user-defined meshsize function $h_u = u(\mathbf{x})$ is provided, this additional constraint is also taken into account at this level. The octant size being bounded from below by the user-defined minimal meshsize $h_{\min}$,
octants are thus subdivided until the condition
\begin{equation}
	h_{octant} \leq \max\left( h_{\min}, \min \left( h_c, h_f, h_u, h_{b}\right) \right)
  \label{eq:refinement}
\end{equation}
is met everywhere in the octree.

\begin{figure}[h!]
	\includegraphics[width=0.49\textwidth]{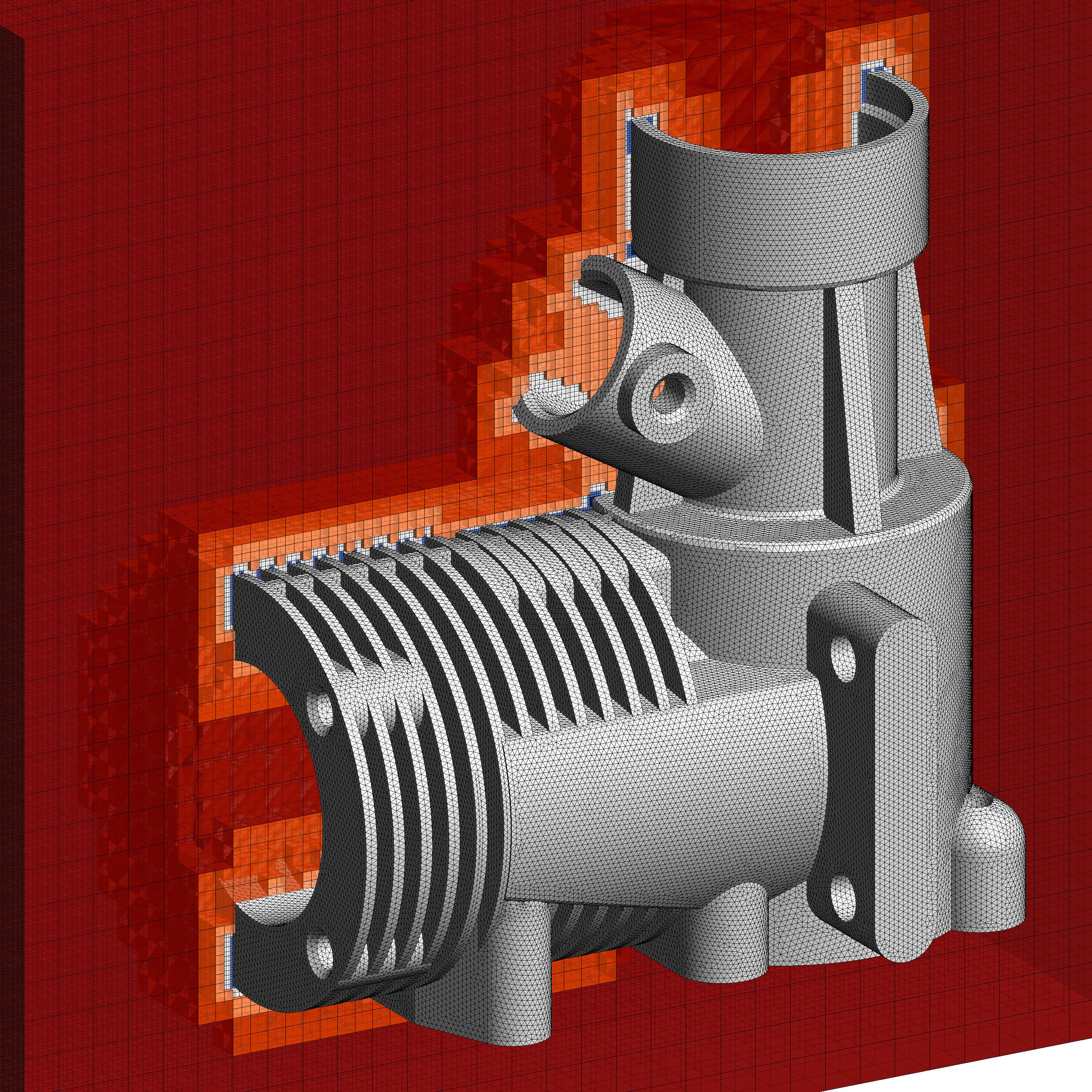}
  \includegraphics[width=0.49\textwidth]{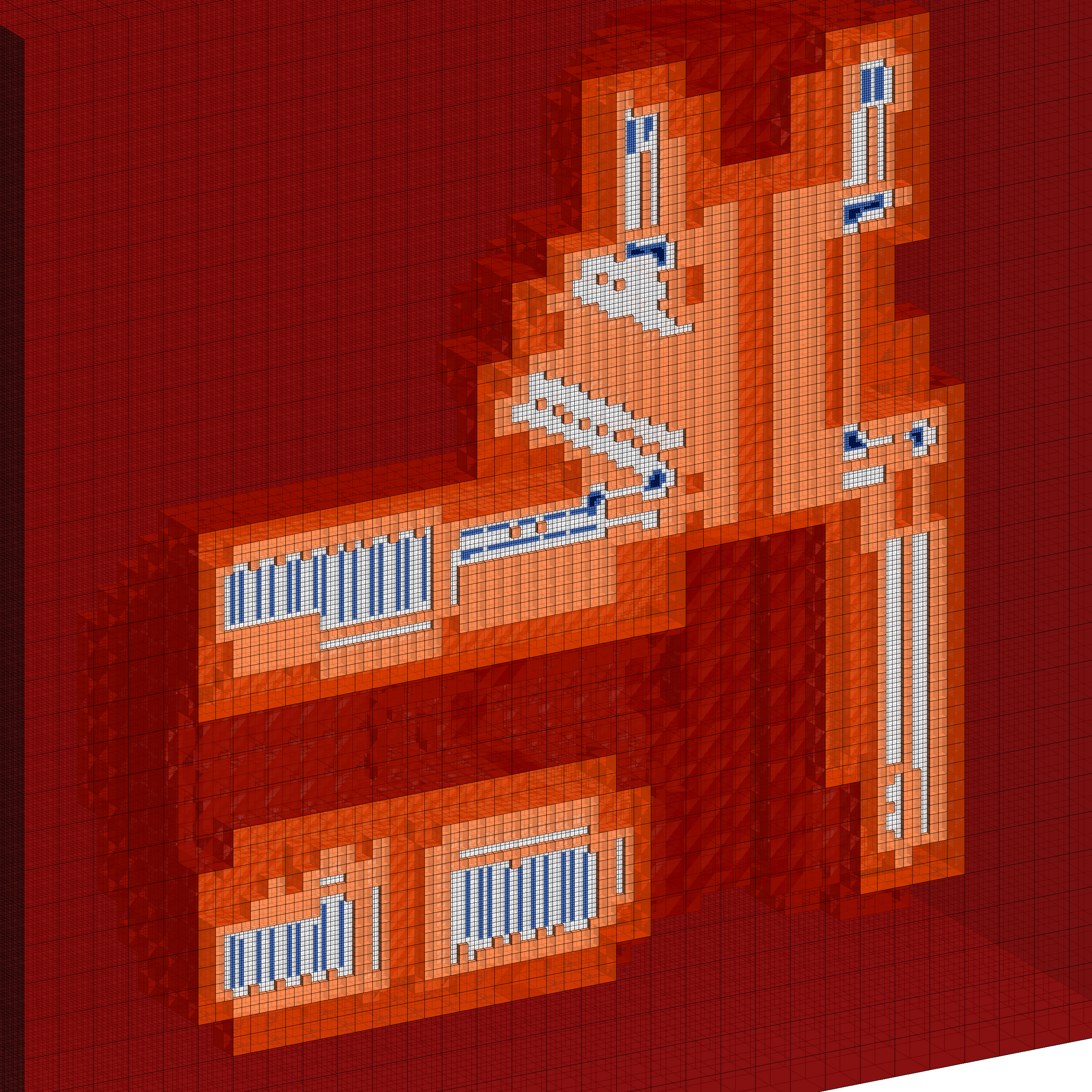}
  \caption{Octree after refinement : the color of each octant signifies its refinement level, from coarse (red) to fine (blue). Octant size is related to the local curvature and feature meshsizes through condition \eqref{eq:refinement}.}
	\label{fig:afterRefinement}
\end{figure}

Once this refinement is completed, the octree is balanced to ensure a maximum 2:1 level ratio between two octants across adjacent faces (Fig. \ref{fig:afterRefinement}), that is, the levels of two octants on each side of a face should not differ by more than one unity. This balancing is necessary to obtain suitable stencils for finite difference computations during the gradient limiting step, see Section 2.4.
It is performed by \textsc{p4est}.
Newly created octants having an intersection with the surface mesh are assigned the meshsize determined by \eqref{eq:refinement}, otherwise the meshsize is set to the bulk size $h_b$ (Fig.~\ref{fig:initialAndSmoothing1}).
At this stage, the meshsize field features large variations,
and a smoothing step is required to end up with a meshsize field suitable for high-quality mesh generation.

\subsection{Regularization}
Large meshsize gradients may be the cause of low quality finite
element solutions. In solid mechanics e.g., they may be the cause of
excessive values of strains when a very small element is adjacent to a
large one. One of the goals of our approach is to ensure that two
adjacent edges in the mesh have their length ratio that are
controlled by a user-defined factor $\alpha$, called the
\emph{regularization parameter}. We now show that this amounts to
limiting the size field gradient by a factor $\alpha-1$ :
\begin{equation}
  \vert \nabla h(\mathbf{x})\vert \leq \alpha - 1.
  \label{smoothing}
\end{equation}

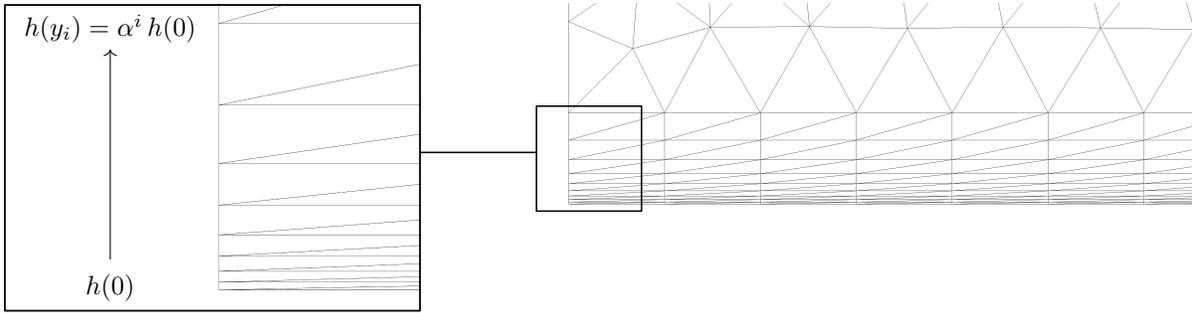
\begin{figure}[h!]
	\centering
		\begin{tikzpicture}
    	\node[anchor=south west,inner sep=0] (image) at (0,0)
    	{\includegraphics[width=\textwidth]{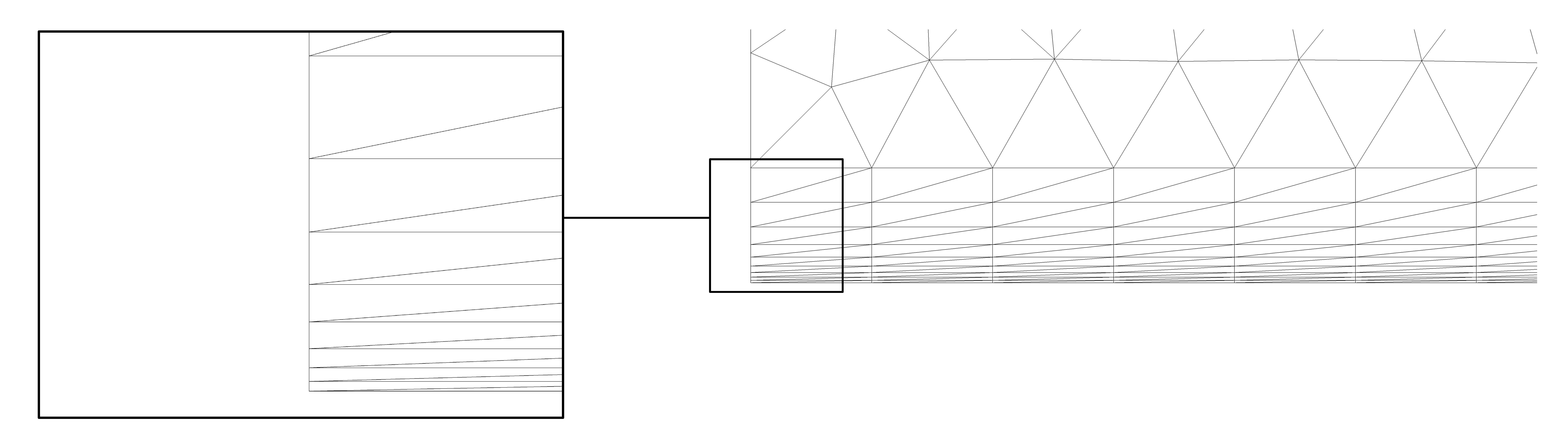}};
    	\begin{scope}[x={(image.south east)},y={(image.north west)}]
				\node at (0.11,0.12) {$h(0)$};
				\node at (0.11,0.86) {$h(y_i) = \alpha^i\,h(0)$};
        \draw[->] (0.11, 0.2) -- (0.11, 0.8);
    	\end{scope}
	\end{tikzpicture}
	\caption{Boundary layer mesh with wall size $h(0) = 0.005$ and gradation $\alpha= 1.4$. The size function is linear with respect to the distance to the wall and grows from $h(0)$ to $\alpha^i\,h(0)$.}
	\label{fig:bLayer}
\end{figure}

This condition may look odd, as one would expect a linear meshsize field to yield a linear progression of ratio $\alpha$ of element sizes, instead of a geometric progression. Let us consider a standard boundary layer mesh (Fig. \ref{fig:bLayer}) defined as follows : a wall size $h(0)$ defines thickness of the first element layer $\delta_0 = h(0)$ and a ratio $\alpha$ defines the geometric progression of element sizes. The $i$-th layer of elements thus has a size $\delta_i = \alpha^i h(0)$. Let us show that a geometric progression of element sizes actually corresponds to an affine size field.
Let $y$ be the vertical distance to the wall.
The coordinate $y_i$ corresponding to the bottom of the $i$-th layer is given by
\begin{equation*}
y_i = \sum_{j=0}^{i-1} \delta_j = h(0) \, \biggr(1 + \alpha + \alpha^2 + \dots + \alpha^{i-1}\biggr) =h(0) \, \frac{\alpha^i - 1}{\alpha-1}.
\end{equation*}
We can thus compute:
$$\alpha^i = 1 + \frac{y_i}{h(0)} (\alpha-1).$$

Since $\delta_i = \alpha^i h(0)$, we have
\begin{equation*}
h(y_i) = \delta_i = h(0) + y_i (\alpha-1)
\end{equation*}
and the mesh size $h(y)$ is affine with respect to the distance to the
wall, so its gradient $\nabla h$ is constant and equal to $\alpha-1$.
A similar computation by Chen et al. [ref] leads to the condition
$\vert \nabla h(\mathbf{x})\vert \leq \ln \alpha$ : both conditions
are very close for typical values of $\alpha \in [1, 1.4]$
\footnote{We have $\ln(1+x) \simeq x$}. As pointed
out by [ref], this one-dimensional analysis is not sufficient to
ensure a global gradation in the mesh. Indeed, when storing the size
function in a background mesh, constraining the size gradient along
edges only ensures the expected gradation along the element edges but
allows for size variation in the interior of the elements in 2D and
3D. To constrain the size gradient on both the edges and the interior,
one has to iterate over the edges of the background mesh and correct
locally the size information stored at the nodes
\cite{chen2017automaticSurf, chen2017automaticVol}. This results in a
size function that satisfies \eqref{smoothing} in a mesh-dependent
fashion \cite{persson2006mesh}. This is not necessary in our
methodology, since the size function is stored in a balanced octree :
condition \eqref{smoothing} is satisfied everywhere in the volume by
iterating over the octants until convergence. This ensures the
gradation is respected in the interior of the elements.

\begin{figure}[h!]
  \centering
	\begin{tikzpicture}
		\pgfmathsetmacro{\cubex}{2}
		\pgfmathsetmacro{\cubey}{2}
		\pgfmathsetmacro{\cubez}{2}
		\draw (0,0,0) -- ++(-\cubex,0,0) -- ++(0,-\cubey,0) -- ++(\cubex,0,0) -- cycle;
		\draw (0,0,0) -- ++(0,0,-\cubez) -- ++(0,-\cubey,0) -- ++(0,0,\cubez) -- cycle;
		\draw (0,0,0) -- ++(-\cubex,0,0) -- ++(0,0,-\cubez) -- ++(\cubex,0,0) -- cycle;
		\draw (-2,0,-2) -- (-2,-2,-2);
		\draw (-2,-2,-2) -- (4,-2,-2);
		\draw (-2,-2,-2) -- (-2,-2,0);

		\draw (2,0,0) -- ++(-\cubex,0,0) -- ++(0,-\cubey,0) -- ++(\cubex,0,0) -- cycle;
		\draw (2,0,0) -- ++(0,0,-\cubez) -- ++(0,-\cubey,0) -- ++(0,0,\cubez) -- cycle;
		\draw (2,0,0) -- ++(-\cubex,0,0) -- ++(0,0,-\cubez) -- ++(\cubex,0,0) -- cycle;

		\draw (4,0,0) -- ++(-\cubex,0,0) -- ++(0,-\cubey,0) -- ++(\cubex,0,0) -- cycle;
		\draw (4,0,0) -- ++(0,0,-\cubez) -- ++(0,-\cubey,0) -- ++(0,0,\cubez) -- cycle;
		\draw (4,0,0) -- ++(-\cubex,0,0) -- ++(0,0,-\cubez) -- ++(\cubex,0,0) -- cycle;

		\node at  (-1,-1,-1) {$h_{i-1}$};
		\node at  (1,-1,-1) {$h_{i}$};
		\node at (3,-1,-1) {$h_{i+1}$};
    \node at (1.3,-2.5,0){FFF};
	\end{tikzpicture}
  \hfill
	\begin{tikzpicture}
		\pgfmathsetmacro{\cubex}{2}
		\pgfmathsetmacro{\cubey}{2}
		\pgfmathsetmacro{\cubez}{2}
		\draw (0,0,0) -- ++(-\cubex,0,0) -- ++(0,-\cubey,0) -- ++(\cubex,0,0) -- cycle;
		\draw (0,0,0) -- ++(0,0,-\cubez) -- ++(0,-\cubey,0) -- ++(0,0,\cubez) -- cycle;
		\draw (0,0,0) -- ++(-\cubex,0,0) -- ++(0,0,-\cubez) -- ++(\cubex,0,0) -- cycle;
		\draw (-2,0,-2) -- (-2,-2,-2);
		\draw (-2,-2,-2) -- (4,-2,-2);
		\draw (-2,-2,-2) -- (-2,-2,0);

		\draw (2,0,0) -- ++(-\cubex,0,0) -- ++(0,-\cubey,0) -- ++(\cubex,0,0) -- cycle;
		\draw (2,0,0) -- ++(0,0,-\cubez) -- ++(0,-\cubey,0) -- ++(0,0,\cubez) -- cycle;
		\draw (2,0,0) -- ++(-\cubex,0,0) -- ++(0,0,-\cubez) -- ++(\cubex,0,0) -- cycle;

		\draw (2,0,0) -- (4,0,0);
		\draw (2,0,-2) -- (4,0,-2);
		\draw (2,-2,0) -- (4,-2,0);
		\draw (4,0,-2) -- (4,0,0);
		\draw (4,-2,0) -- (4,0,0);
		\draw (4,-2,0) -- (4,-2,-2);
		\draw (4,0,-2) -- (4,-2,-2);
		\draw (4,0,-1) -- (4,-2,-1);
		\draw[densely dotted] (4,-1,0) -- (4,-1,-2);
		\draw[densely dotted] (2,0,-1) -- (4,0,-1);
		\draw[densely dotted] (2,-2,-1) -- (4,-2,-1);
		\draw[densely dotted] (3,-2,0) -- (3,-2,-2);
		\draw[densely dotted] (3,0,-2) -- (3,-2,-2);
		\draw[densely dotted] (3,0,-1) -- (3,-2,-1);
		\draw[densely dotted] (3,-1,0) -- (3,-1,-2);
		\draw[densely dotted] (2,-1,0) -- (4,-1,0);
		\draw[densely dotted] (2,0,-1) -- (2,-2,-1);
		\draw[densely dotted] (2,-1,-1) -- (4,-1,-1);
		\draw[densely dotted] (2,-1,-2) -- (4,-1,-2);
		\draw[densely dotted] (2,-1,0) -- (2,-1,-2);
		\draw[densely dotted] (3,0,0) -- (3,-2,0);
		\draw[densely dotted] (3,0,0) -- (3,0,-2);

		\node at (-1,-1,-1) {$h_{i-1}$};
		\node at (1,-1,-1) {$h_{i}$};

		\node at (2.5,-0.5,-0.5) {$h_{1}$};
		\node at (2.5,-1.5,-0.5) {$h_{2}$};
		\node at (2.5,-1.5,-1.5) {$h_{3}$};
		\node at (2.5,-0.5,-1.5) {$h_{4}$};

    \node at (1.3,-2.5,0){FFH};
	\end{tikzpicture}
  \begin{tikzpicture}
		\pgfmathsetmacro{\cubex}{2}
		\pgfmathsetmacro{\cubey}{2}
		\pgfmathsetmacro{\cubez}{2}
    \draw (-2,0,0) -- (0,0,0);
    \draw (-2,0,0) -- (-2,0,-2);
    \draw (-2,0,0) -- (-2,-2,0);
    \draw (-2,0,-2) -- (-2,-2,-2);
    \draw (-2,-2,0) -- (-2,-2,-2);
    \draw (-2,-2,-2) -- (0,-2,-2);
		\draw (-2,0,-2) -- (0,0,-2);
		\draw (-2,-2,0) -- (0,-2,0);
		\draw (0,0,-2) -- (0,0,0);
		\draw (0,-2,0) -- (0,0,0);
		\draw (0,-2,0) -- (0,-2,-2);
		\draw (0,0,-2) -- (0,-2,-2);
		\draw (0,-2,-2) -- (2,-2,-2);
		\draw[densely dotted] (0,0,-1) -- (0,-2,-1);
		\draw[densely dotted] (0,-1,0) -- (0,-1,-2);
		\draw[densely dotted] (-2,0,-1) -- (0,0,-1);
		\draw[densely dotted] (-2,-2,-1) -- (0,-2,-1);
		\draw[densely dotted] (-1,-2,0) -- (-1,-2,-2);
		\draw[densely dotted] (-1,0,-2) -- (-1,-2,-2);
		\draw[densely dotted] (-1,0,-1) -- (-1,-2,-1);
		\draw[densely dotted] (-1,-1,0) -- (-1,-1,-2);
		\draw[densely dotted] (-2,-1,0) -- (0,-1,0);
		\draw[densely dotted] (-2,0,-1) -- (-2,-2,-1);
		\draw[densely dotted] (-2,-1,-1) -- (0,-1,-1);
		\draw[densely dotted] (-2,-1,-2) -- (0,-1,-2);
		\draw[densely dotted] (-2,-1,0) -- (-2,-1,-2);
		\draw[densely dotted] (-1,0,0) -- (-1,-2,0);
		\draw[densely dotted] (-1,0,0) -- (-1,0,-2);
    \node at (-0.5,-0.5,-0.5) {$h_{1}$};
		\node at (-0.5,-1.5,-0.5) {$h_{2}$};
		\node at (-0.5,-1.5,-1.5) {$h_{3}$};
		\node at (-0.5,-0.5,-1.5) {$h_{4}$};

		\draw (2,0,0) -- ++(-\cubex,0,0) -- ++(0,-\cubey,0) -- ++(\cubex,0,0) -- cycle;
		\draw (2,0,0) -- ++(0,0,-\cubez) -- ++(0,-\cubey,0) -- ++(0,0,\cubez) -- cycle;
		\draw (2,0,0) -- ++(-\cubex,0,0) -- ++(0,0,-\cubez) -- ++(\cubex,0,0) -- cycle;

		\draw (2,0,0) -- (4,0,0);
		\draw (2,0,-2) -- (4,0,-2);
		\draw (2,-2,0) -- (4,-2,0);
		\draw (4,0,-2) -- (4,0,0);
		\draw (4,-2,0) -- (4,0,0);
		\draw (4,-2,0) -- (4,-2,-2);
		\draw (4,0,-2) -- (4,-2,-2);
		\draw (4,0,-1) -- (4,-2,-1);
		\draw[densely dotted] (4,-1,0) -- (4,-1,-2);
		\draw[densely dotted] (2,0,-1) -- (4,0,-1);
		\draw[densely dotted] (2,-2,-1) -- (4,-2,-1);
		\draw[densely dotted] (3,-2,0) -- (3,-2,-2);
		\draw[densely dotted] (3,0,-2) -- (3,-2,-2);
		\draw[densely dotted] (3,0,-1) -- (3,-2,-1);
		\draw[densely dotted] (3,-1,0) -- (3,-1,-2);
		\draw[densely dotted] (2,-1,0) -- (4,-1,0);
		\draw[densely dotted] (2,0,-1) -- (2,-2,-1);
		\draw[densely dotted] (2,-1,-1) -- (4,-1,-1);
		\draw[densely dotted] (2,-1,-2) -- (4,-1,-2);
		\draw[densely dotted] (2,-1,0) -- (2,-1,-2);
		\draw[densely dotted] (3,0,0) -- (3,-2,0);
		\draw[densely dotted] (3,0,0) -- (3,0,-2);
		\node at (1,-1,-1) {$h_{i}$};

		\node at (2.5,-0.5,-0.5) {$h_{1}$};
		\node at (2.5,-1.5,-0.5) {$h_{2}$};
		\node at (2.5,-1.5,-1.5) {$h_{3}$};
		\node at (2.5,-0.5,-1.5) {$h_{4}$};

    \node at (1.3,-2.5,0){HFH};
	\end{tikzpicture}
	\caption{Finite difference stencils to compute $\nabla h$ on the background mesh.}
	\label{fig:stencils}
\end{figure}
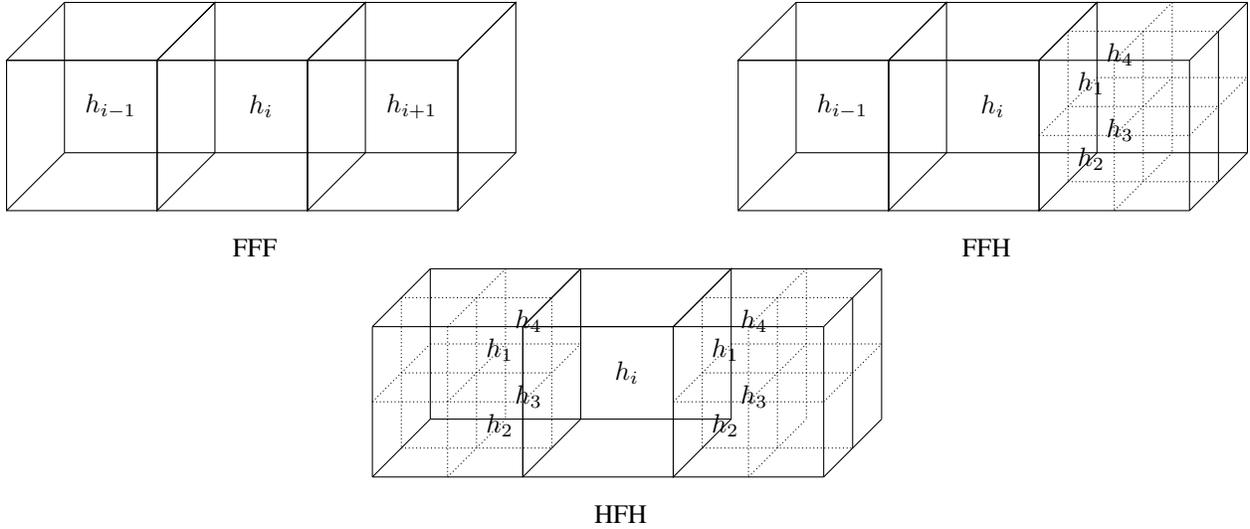

\emph{Gradient computation by finite difference.}
~The regularization parameter $\alpha$ here introduced limits the maximum ratio between
two adjacent edges of the mesh. To ensure that, we iterate over the octants to limit $h$ until
\begin{equation}
	\left\vert \frac{\partial h}{\partial x_i} \right\vert \leq \alpha-1,
\end{equation}
everywhere, with $i = 1, \hdots, 3$. For each octant, we take advantage of the 2:1 balancing
provided by \textsc{p4est} to compute the gradient with a
cell-centered finite difference scheme. Since the octree is balanced,
only three different stencils are necessary to approximate the
derivatives. To use \textsc{p4est} terminology, one side of a face
between two octants is said to be either full (F) or hanging (H),
depending if the octant on that side is a leaf or is itself
subdivided, respectively. The three stencils (Fig. \ref{fig:stencils}) are then FFF, FFH, HFH,
the all-hanging case being divided in multiple all-full stencils. To evaluate the derivatives of the sizing function at the center of the middle octant along, say, the \emph{x}-direction, we use the following Taylor approximation :
\begin{equation}
  \frac{\partial h}{\partial x}\biggr\vert_i = \frac{\bar{h}_{i+1} - h_i}{2\Delta x_{i+1}} + \frac{h_i-\bar{h}_{i-1}}{2\Delta x_{i-1}},
\end{equation}
where $\bar{h}$ denotes the average of the size stored in the four adjacent octants to a hanging face, and $\Delta x$ is the sum of the half-lengths of octants on each side of a face. For the all-full stencil, $\bar{h} = h$ and $\Delta x$ is a constant step equal to the length of an octant, so that this approximation reduces to the usual second-order centered scheme. By taking the average size $\bar{h}$ in the hanging case, we decouple the equations to solve for the discrete approximation of the size gradient, and avoid solving a least-square problem at each octant and at each iteration.

\begin{figure}[h!]
\centering
\begin{subfigure}{.49\textwidth}
  \centering
  \includegraphics[width=\linewidth]{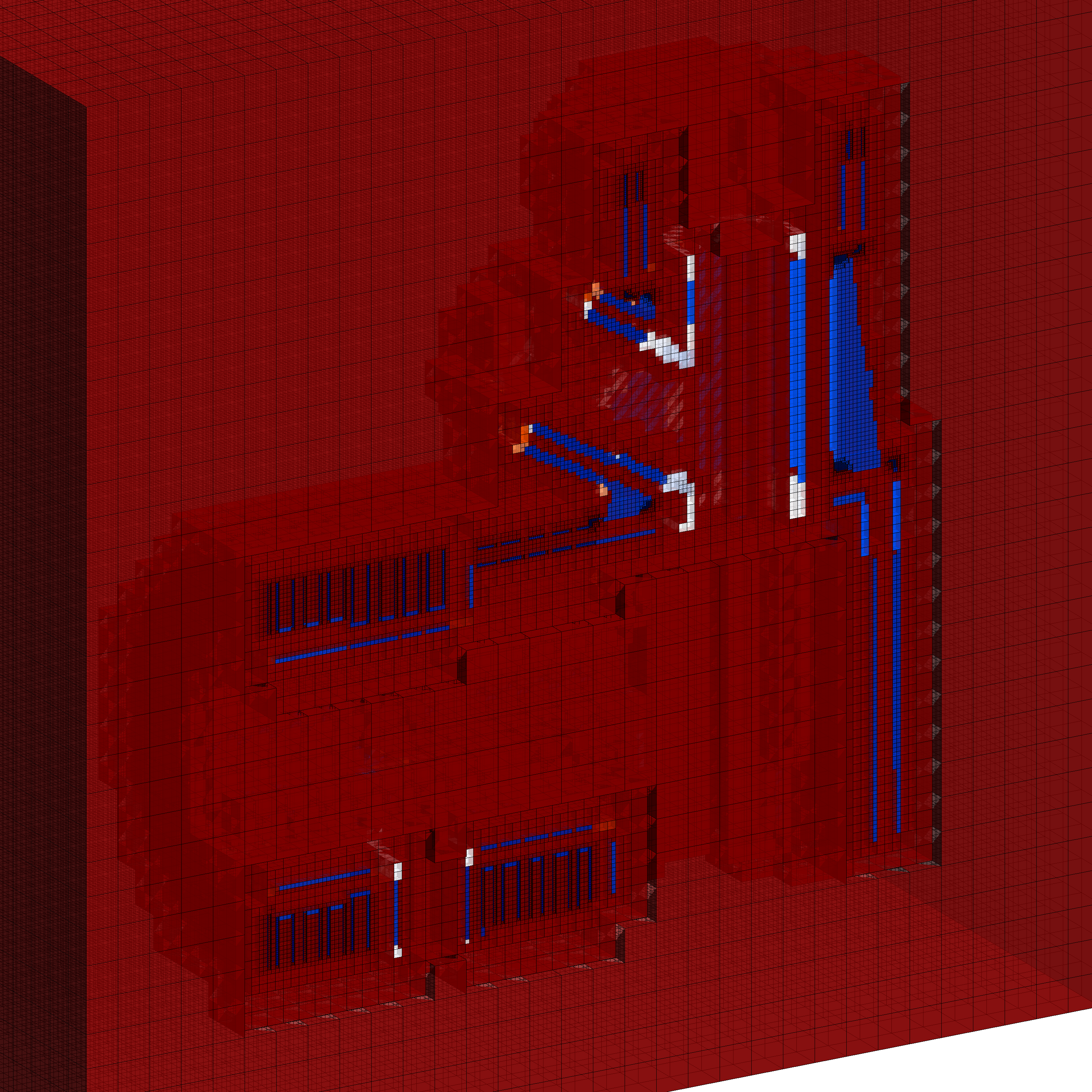}
  \caption{}
  \label{fig:initialAndSmoothing1}
\end{subfigure}
\hfill
\begin{subfigure}{.49\textwidth}
  \centering
  \includegraphics[width=\linewidth]{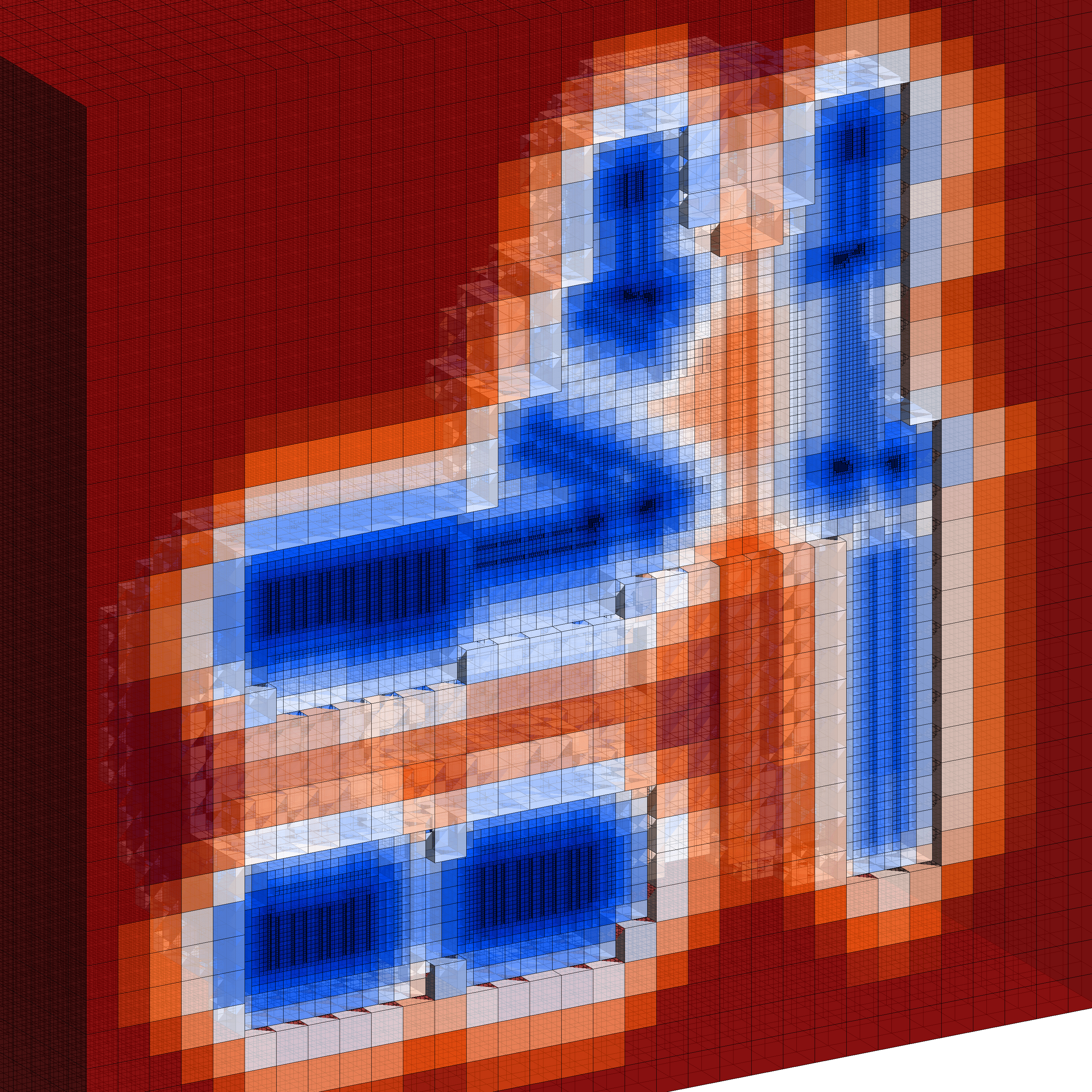}
  \caption{}
  \label{fig:initialAndSmoothing2}
\end{subfigure}

\caption{Limiting the meshsize stored in the octants : (a) initial meshsize computed from curvature and feature size, assigned in octants intersecting the surface mesh (b) meshsize after regularization.}
\label{fig:initialAndSmoothing}
\end{figure}

\emph{Size limitation.} ~The computed gradient is then used to limit the size field. For any two adjacent octants with gradient larger than $\alpha-1$, larger sizes are limited in order to satisfy condition \eqref{smoothing} : let $h_1$ and $h_2$ denote the size stored on the octants, with $h_2 \geq h_1$. The larger size $h_2$ is corrected using the linear extrapolation :
\begin{equation}
  h_2 = \min(h_2, \,h_1 + \Delta x (\alpha-1))
\end{equation}
This expression is the discretized version of the steady-state solution to the Hamilton-Jacobi equation proposed by Persson \cite{persson2006mesh} in their continuous formulation of the gradient limiting problem, ensuring size limitation is propagated in the direction of increasing values. Smaller sizes are left unchanged, so that sharp features of the geometry are preserved in the size function and are not overlooked during surface or volume meshing. The smoothing is performed iteratively in the three directions, until the constraint $\vert \nabla h \vert \leq \alpha-1$ is satisfied everywhere in the octree, propagating the small sizes and yielding a continuous size field (Fig.\ref{fig:initialAndSmoothing2}).

\subsection{Size query in the octree}
Size queries are performed by Gmsh to evaluate the size at different locations on parametrized curves, surfaces and volumes during the meshing process. The implementation of queries routine in the octree is provided by \textsc{p4est}. Since both the size and its gradient are known in each octant, the size at the query point $\mathbf{x}$ is approximated by a first-order Taylor expansion, that is : $$h(\mathbf{x}) = h_i + \nabla h \cdot (\mathbf{x} - \mathbf{x_c}),$$ where $h_i$ and $\mathbf{x}_c$ denote the size at and the coordinate of the center of the octant, respectively.

\section{Results}

To illustrate our methodology, our algorithm has been applied on a
large variety CAD models. The models were found at different
locations: GrabCAD (https://grabcad.com/), the ABC Dataset
model library \cite{Koch_2019_CVPR} and Gmsh's benchmarks suite.
For each of the following test cases, we start by tessellating the CAD
model using a uniform mesh size over the whole geometry. The
resulting triangular mesh is the input of our algorithm. We then
compute (i) surface curvatures and (ii) the approximate medial
axis. Those two informations are processesed to build a size field as
described in the previous sections of the paper. The octree is stored
on disk in the native \textsc{p4est} format and is loaded
 as a background field in Gmsh \footnote{{\tt gmsh
  foo.geo -bgm foo.p4est}}.

This sections aims at showing that our approach is both efficient and
accurate. As far as efficiency is concerned, we distinguish the CPU
time dedicated to queries in the octree and the CPU time devoted to
the construction of the octree. We then compare the sum of those two
to the remaining time for meshing. Note here that the last evolution of
Gmsh's 3D mesh generator is extremely efficient\cite{marot2019one}. We
show that the new size field does not change orders of magnitude in
the meshing time. In other word, seconds remain seconds and do not
become minutes. With regard to accuracy,
we check wether critical parameters of our approach are taken ino
account accurately. Discrete gradients of mesh sizes are compared with their
limit value $\alpha-1$. The number of elements in small gaps is
checked to be close to $n_g$. Finally,  the efficiency index of size
fields, which measures the discrepancy between adimensional edge
lengths and their ideal values $1$, is presented.

All computations are run on a laptop with Intel Core i7 8750h CPU (2.2
GHz) and 16Gb memory. Execution times for two selected test cases can
be found in the table hereafter.

The parameters chosen for the size field constructionss are :
\begin{itemize}
  \item the bulk size $h_b$ is set to $L/20$, where $L$ denotes the maximum dimension of the axis-aligned bounding box of the model;
  \item the minimum size $h_{\min}$ is set to $L/1000$;
  \item the node density $n_d$ is set to $20$;
  \item the number of layers in geometric features $n_g$ is set to $4$;
  \item the gradation $\alpha$ is set to $1.1$.
\end{itemize}

\subsection{Surface and volume meshing}
\begin{figure}
	\includegraphics[width=0.95\textwidth]{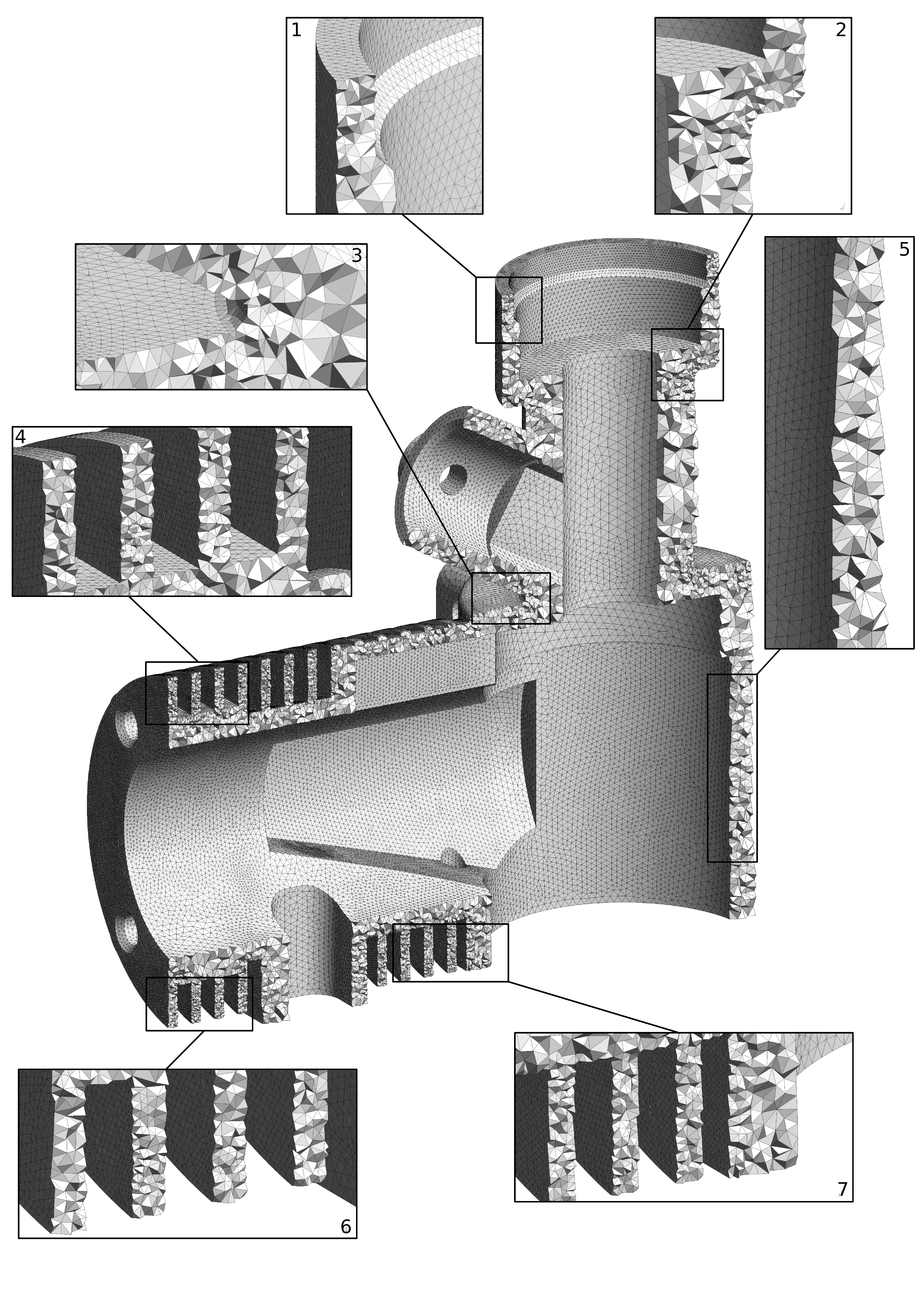}
	\caption{Final volume mesh of our worked-out example : the size field is computed with $n_d = 20$ nodes on the local osculating circles and $n_g = 4$ layers of elements in all geometric features. The mesh contains 583,776 nodes and 2,716,170 tetrahedra.}
	\label{fig:blockVolumeMesh}
\end{figure}

\emph{Running example.}~We first present the resulting meshes of the
engine block test case. The geometry, while being relatively simple,
has all features that should be captured by our approach: it presents entities of variable
radii of curvature each endowed with a specific curvature meshsize. The
CAD model is mostly an assemblage of thin components such as thin cylinders and narrow fins surrounding the main cylinder: feature meshsize is thus the dominant meshing criterion
overall. Surface meshes are presented in Section 2, and show the
impact of the feature meshsize (Fig. \ref{fig:withMedialAxis}). On the left part of Figure \ref{fig:withMedialAxis}, the meshsize field is based only on the curvature meshsize $h_c$: we see that only one or two elements are generated in the narrow fins, whose meshsize is mostly defined by the radius
of curvature of the inner cylinder. On the right part of Figure
\ref{fig:withMedialAxis}, feature meshsize is included in the
meshsize field and we observe $n_g = 4$ layers in the fins.
Figure \ref{fig:blockVolumeMesh} reveals the tetrahedral mesh and
confirms that meshsize specifications are respected inside the
volume by the meshing tool. Meshsize in curved areas, such as
corners, is constrained by curvature rather than geometric features,
since branches of the medial axis were removed in their vicinity. This
results in a smaller meshsize defined by the node density $n_d$
(zooms 2 and 3). Meshsize in the thin cylinders and the fins is
computed to fit 4 layers of elements (zooms 1, 4, 5, 6 and 7).

\begin{figure}[h!]
\centering
  \includegraphics[width=0.49\linewidth]{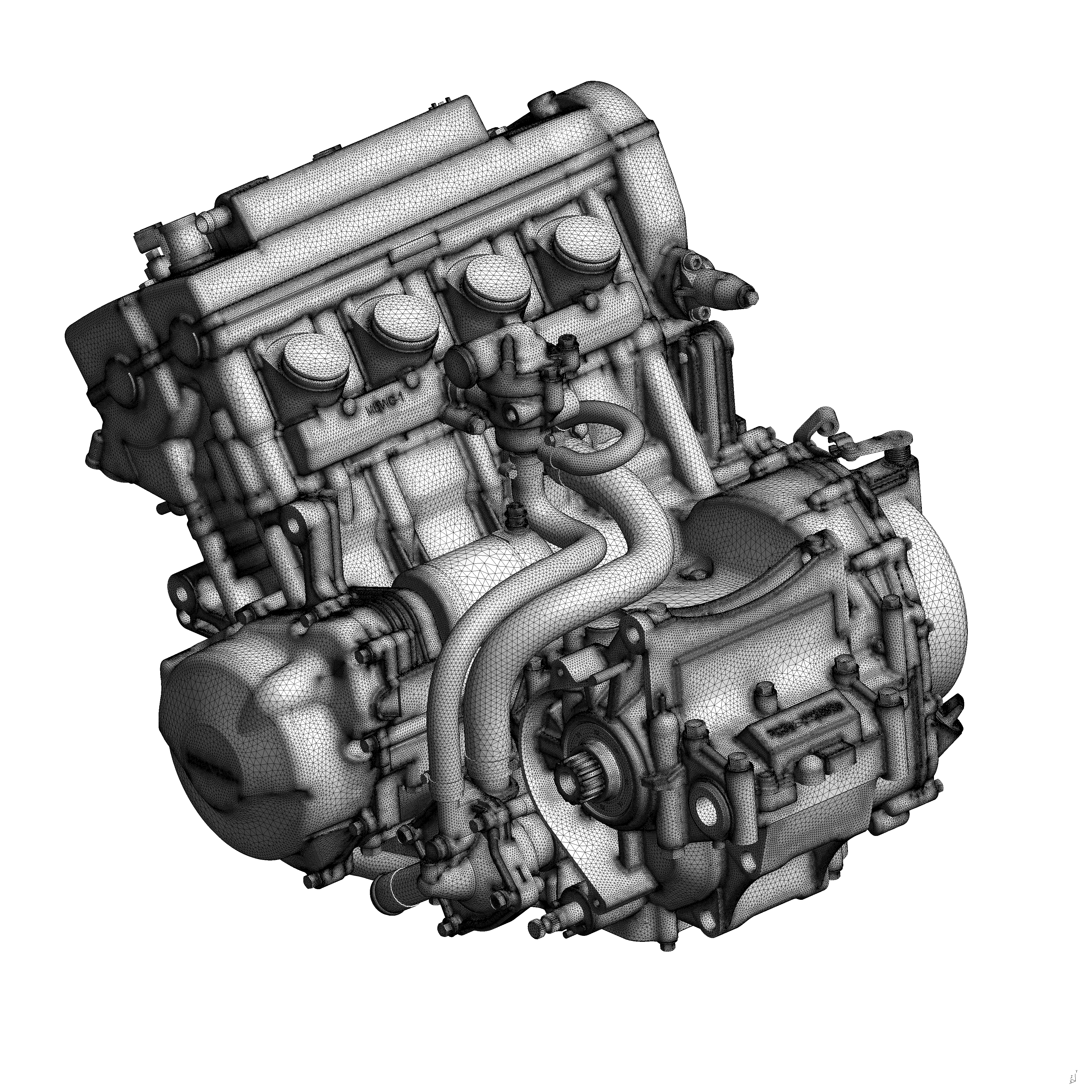}
  \includegraphics[width=0.49\linewidth]{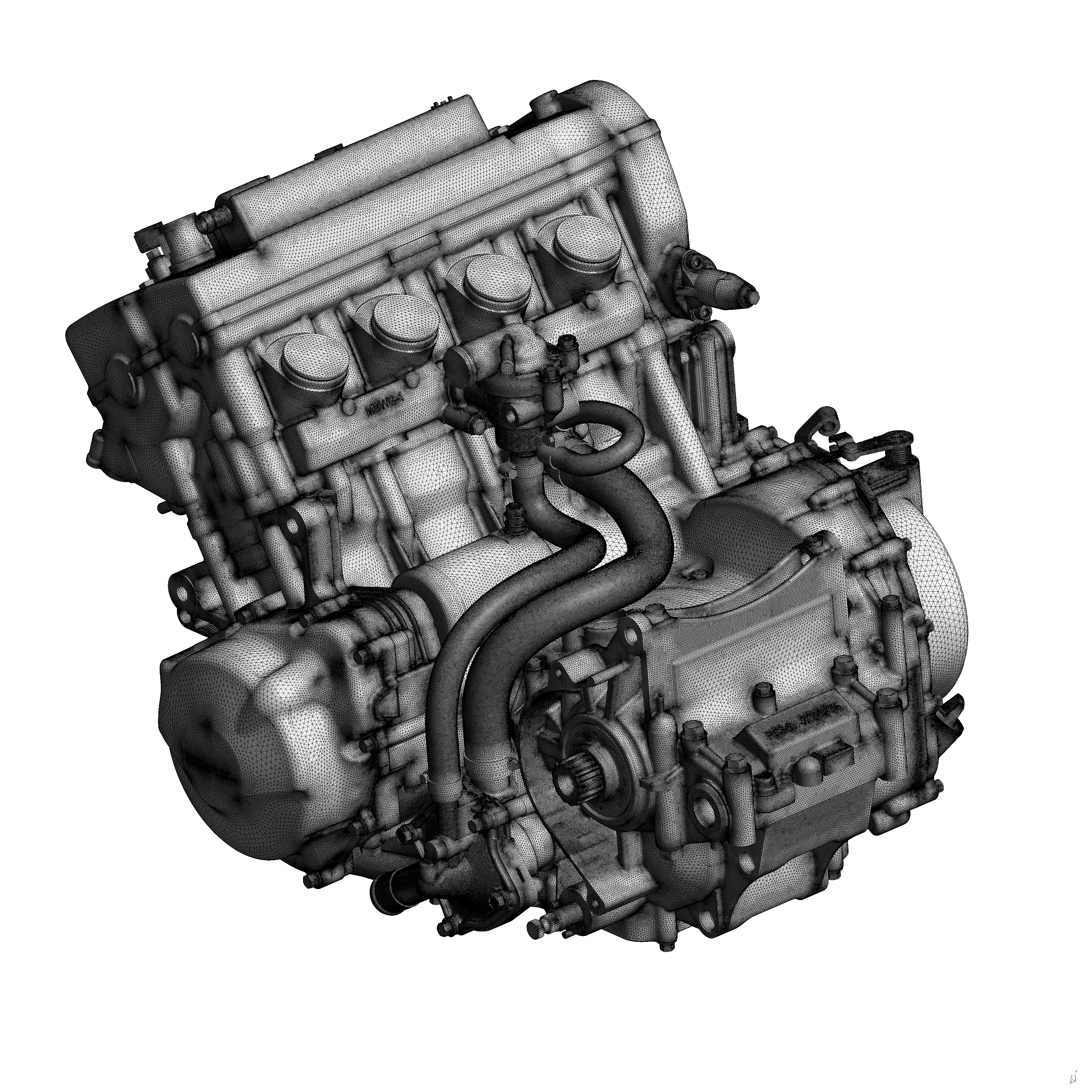}
  \caption{The surface mesh based only on curvature contains 1,521,410 nodes and 3,042,060 triangles (left); the mesh based on both curvature and feature size contains 2,302,630 nodes and 4,605,010 triangles (right).}
  \label{fig:honda_mesh}
\end{figure}

\begin{figure}[h!]
\centering
  \includegraphics[width=0.49\linewidth]{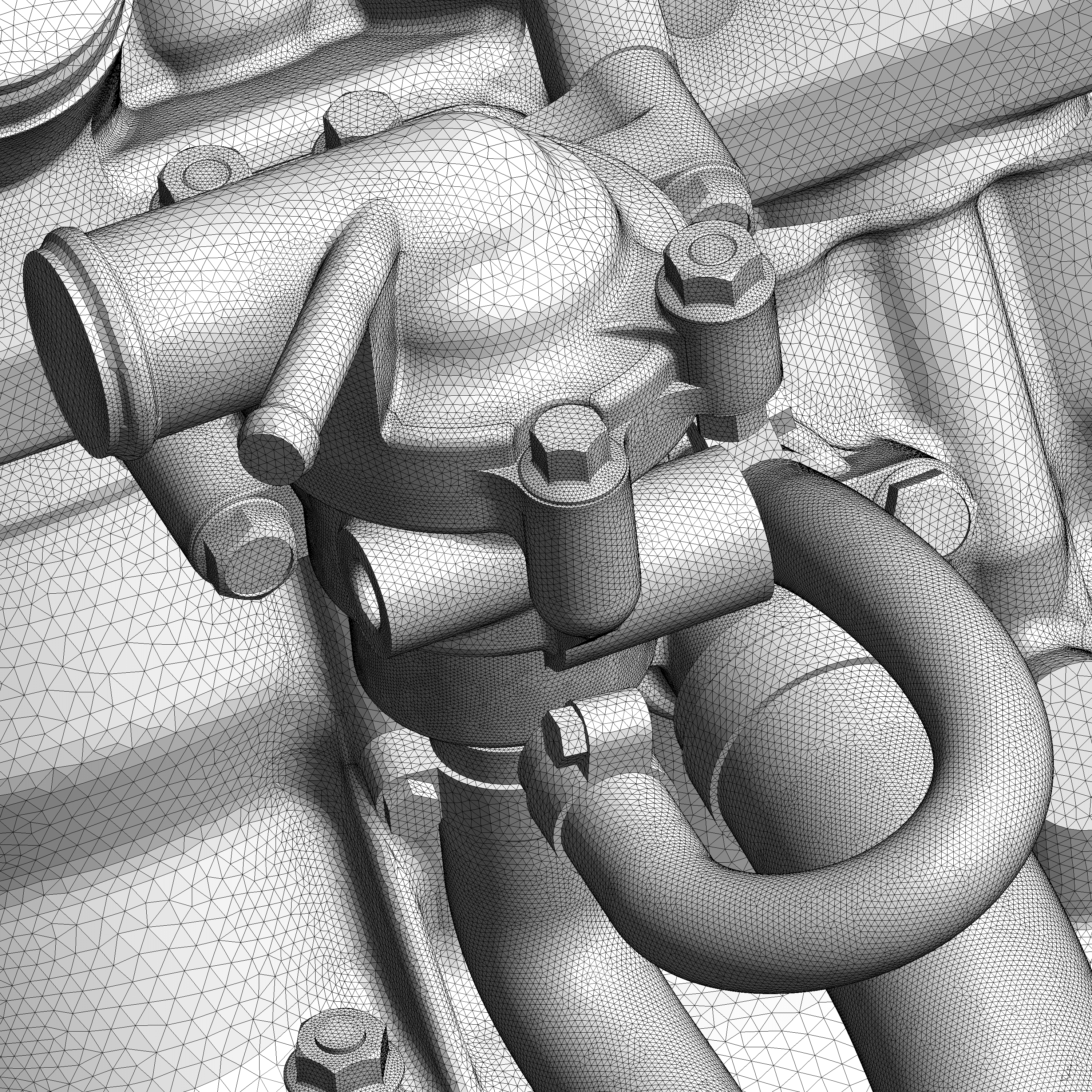}
  \includegraphics[width=0.49\linewidth]{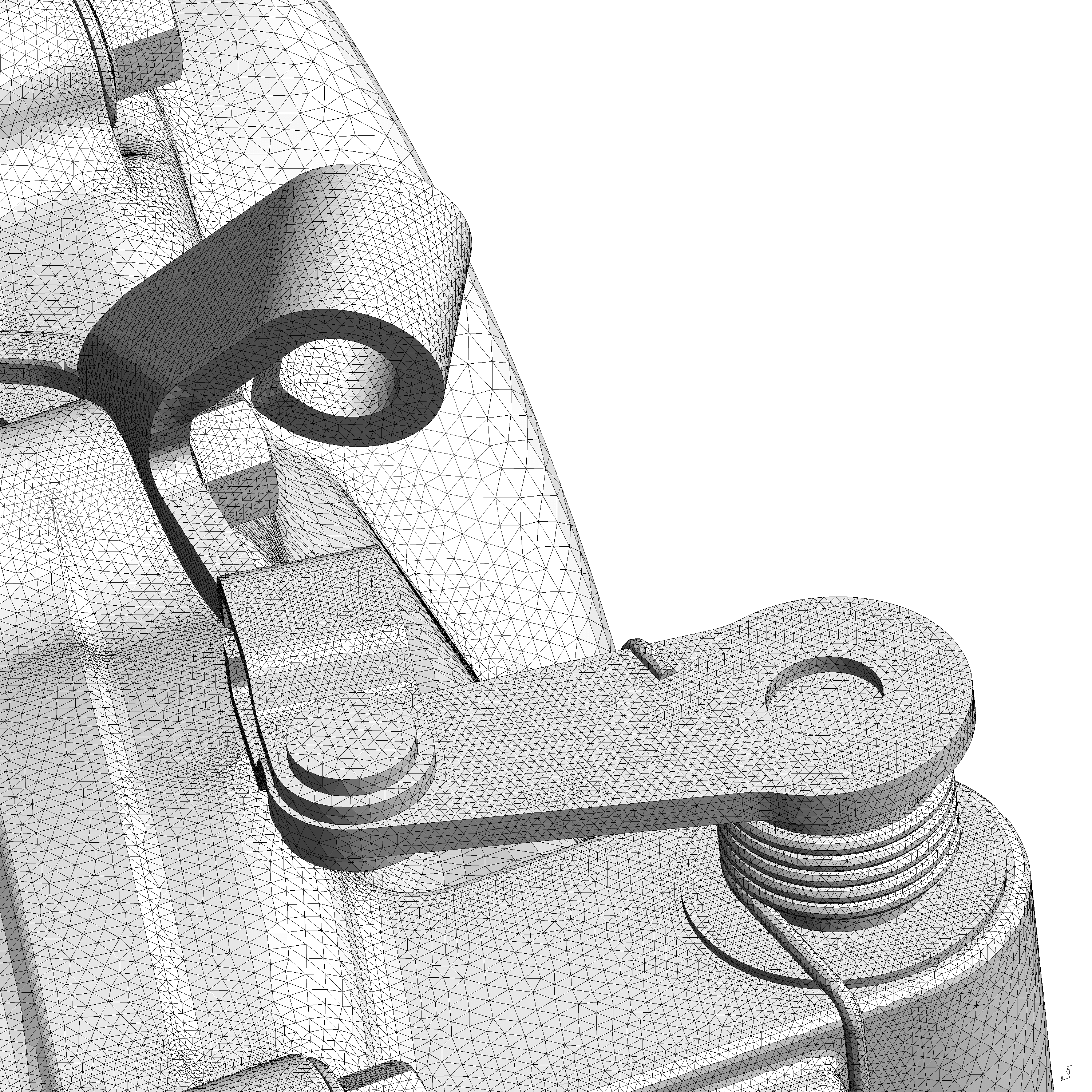}
  \includegraphics[width=0.49\linewidth]{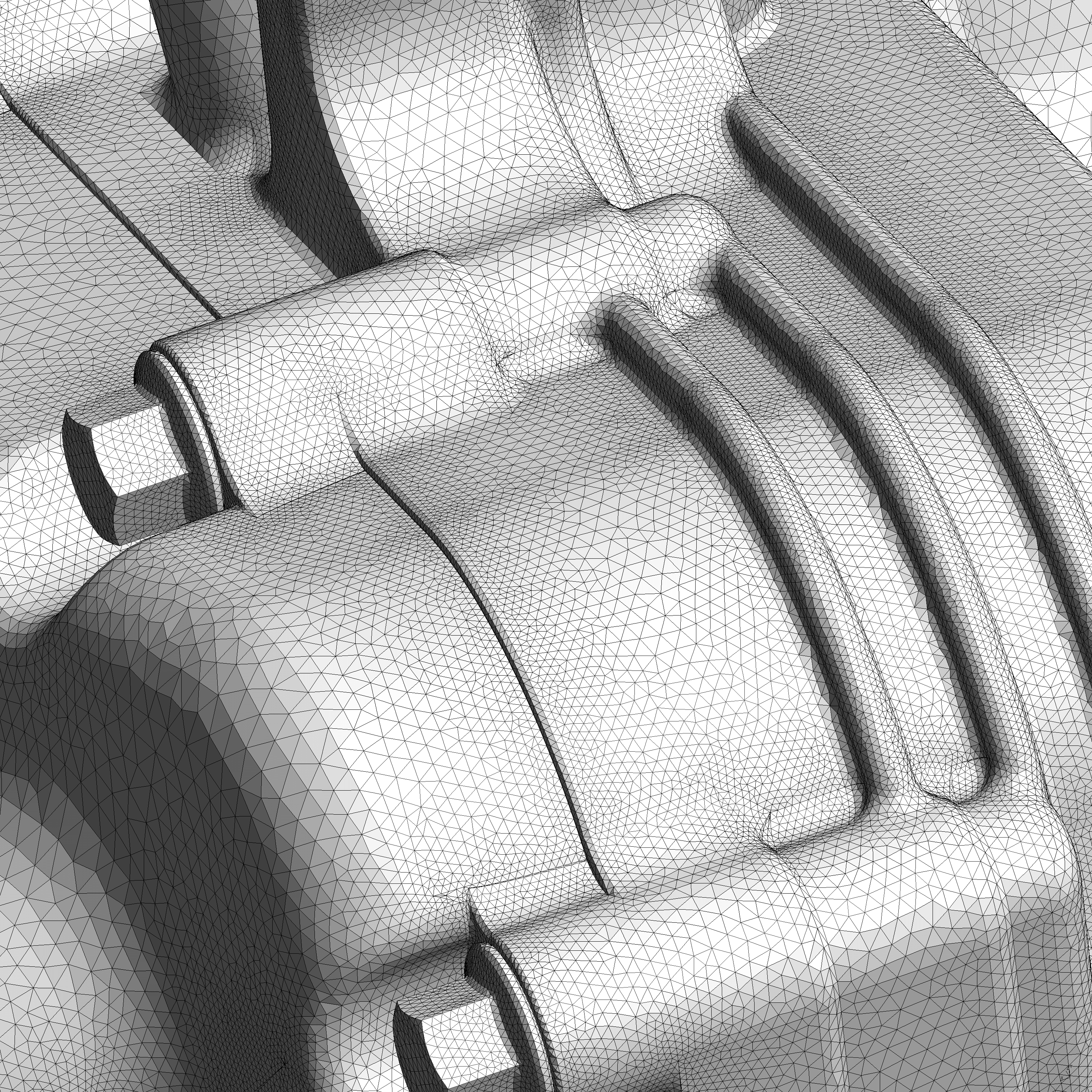}
  \includegraphics[width=0.49\linewidth]{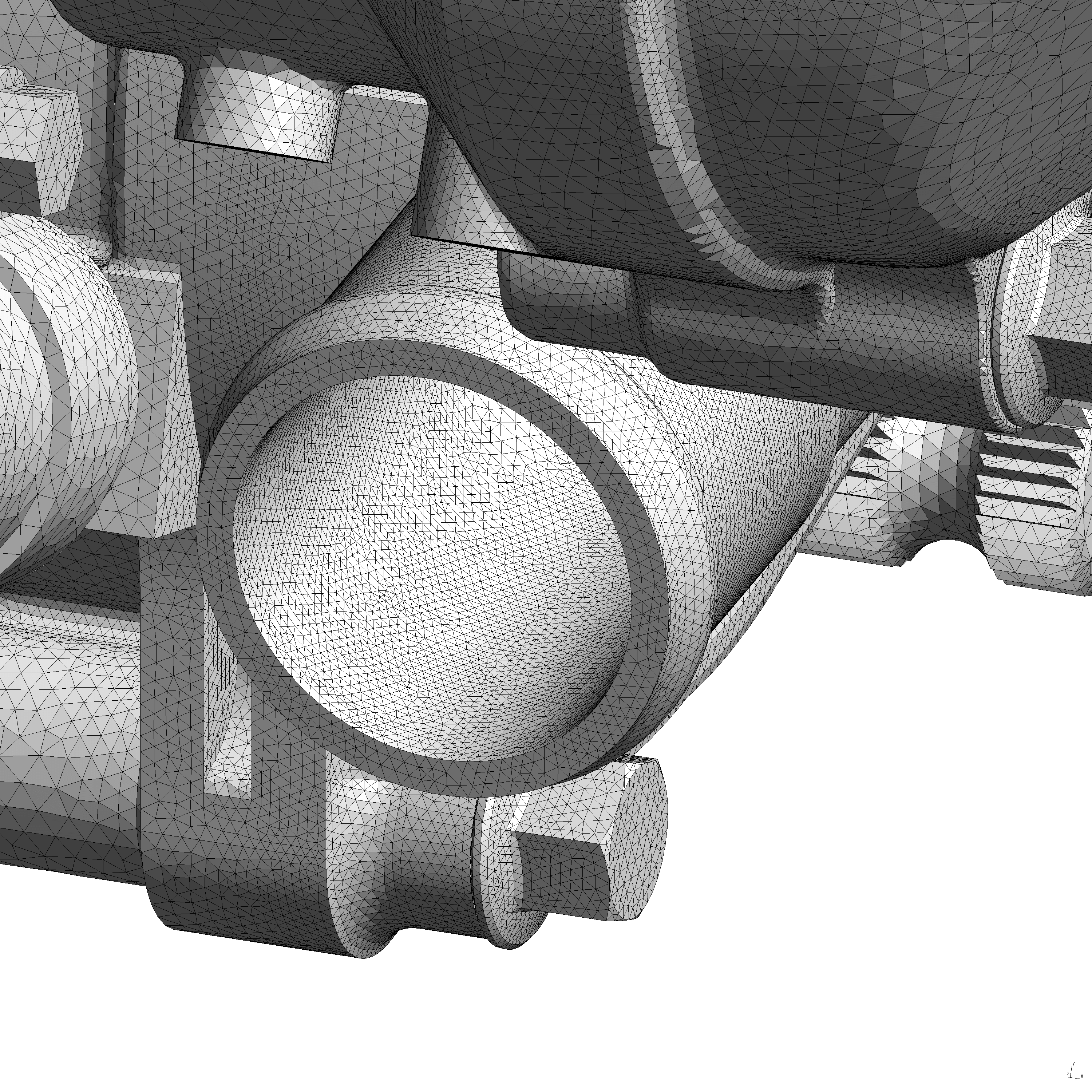}
  \caption{Zooms on selected parts of the Honda engine.}
  \label{fig:honda_mesh_zoom}
\end{figure}

\emph{Sport bike engine.}~Our second test case is the four-cylinder
engine of a Honda{\texttrademark} CBR600F4i sport bike, accessed from
GrabCAD\footnote{https://grabcad.com/library/honda-cbr600-f4i-engine-1}. The
model is composed of thin areas such as pipes, gears and plates, and
of a large number of features overall, making it a interesting test
case for our algorithm. Unfortunately, the volume enclosed is not
"air-tight", and can not be meshed without repairing the CAD model,
hence only surface meshes are presented for this example. The
resulting surface meshes accurately capture areas of higher curvature
(Fig. \ref{fig:honda_mesh}, left) as well as the small features
(Fig. \ref{fig:honda_mesh}, right and Fig. \ref{fig:honda_mesh_zoom}),
demonstrating the robustness of our algorithm on such large CAD
models.


\begin{figure}[h!]
\centering
\begin{subfigure}{.3\textwidth}
  \centering
  \includegraphics[width=\linewidth]{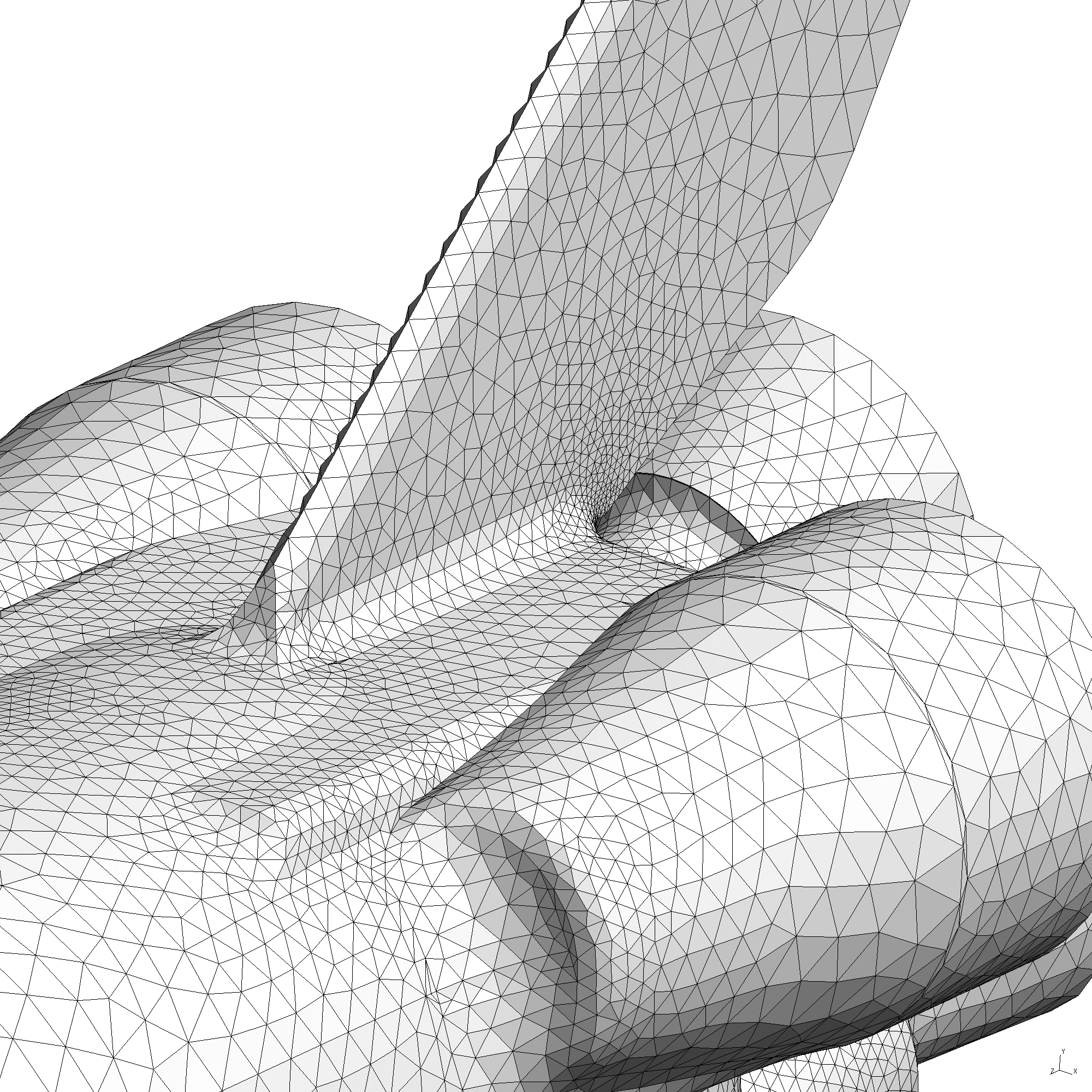}
  \caption{$n_d = 10$}
\end{subfigure}%
\hfill
\begin{subfigure}{.3\textwidth}
  \centering
  \includegraphics[width=\linewidth]{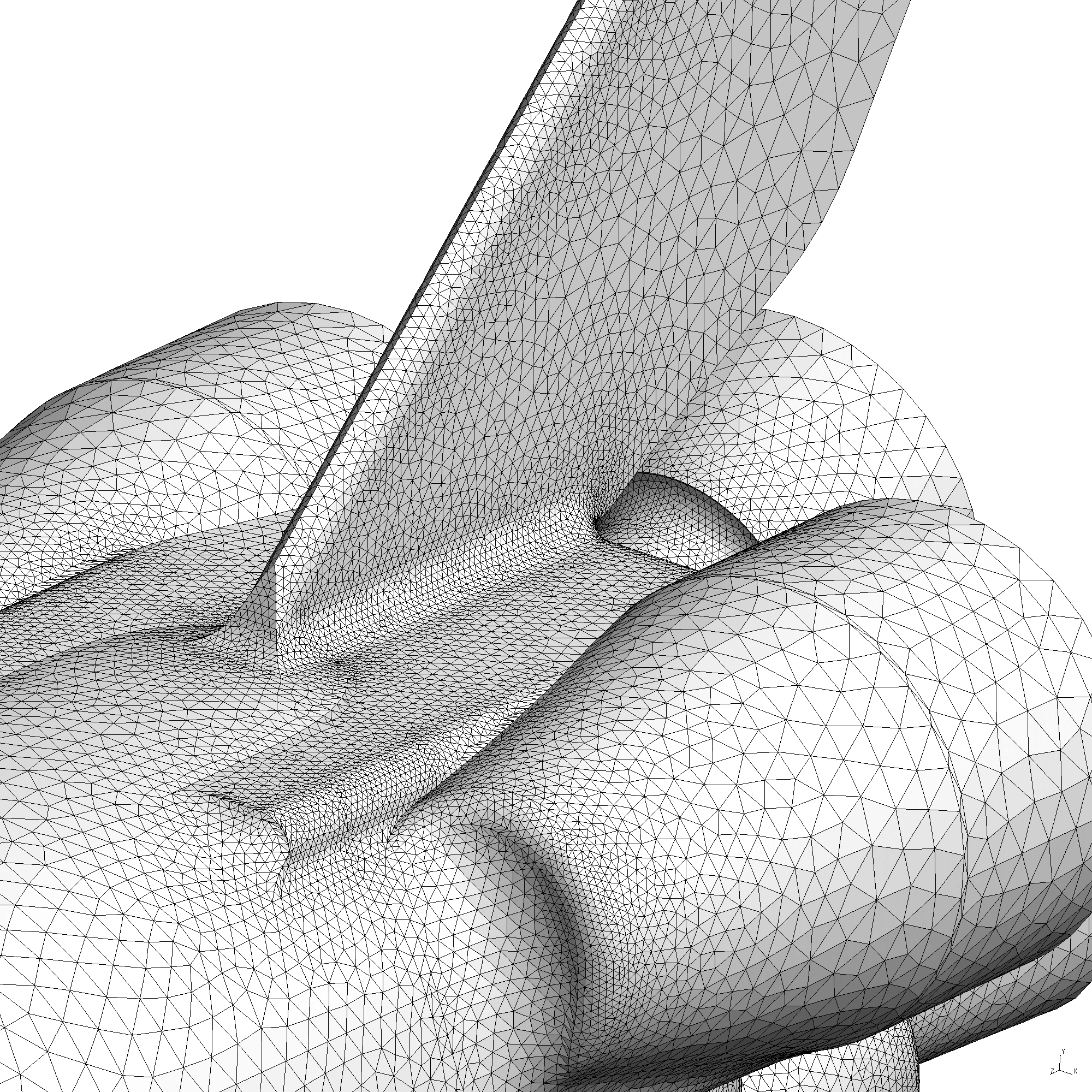}
  \caption{$n_d = 30$}
\end{subfigure}
\centering
\hfill
\begin{subfigure}{.3\textwidth}
  \centering
  \includegraphics[width=\linewidth]{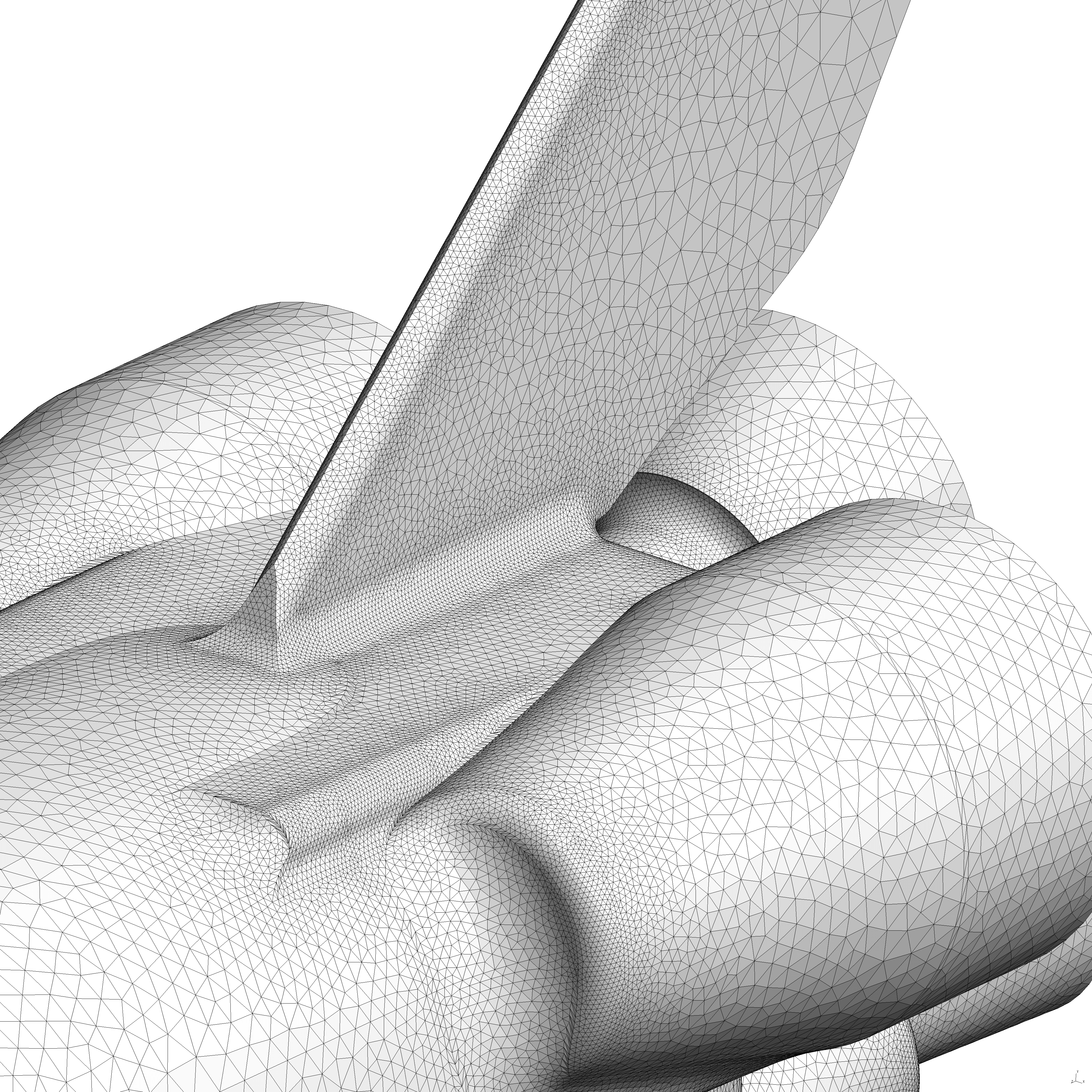}
  \caption{$n_d = 50$}
\end{subfigure}%
\caption{Influence of the node density $n_d$ on the final mesh for a gradation of $\alpha = 1.1$.}
\label{fig:influenceDensity}
\end{figure}

\begin{figure}[h!]
\centering
\begin{subfigure}{.3\textwidth}
  \centering
  \includegraphics[width=\linewidth]{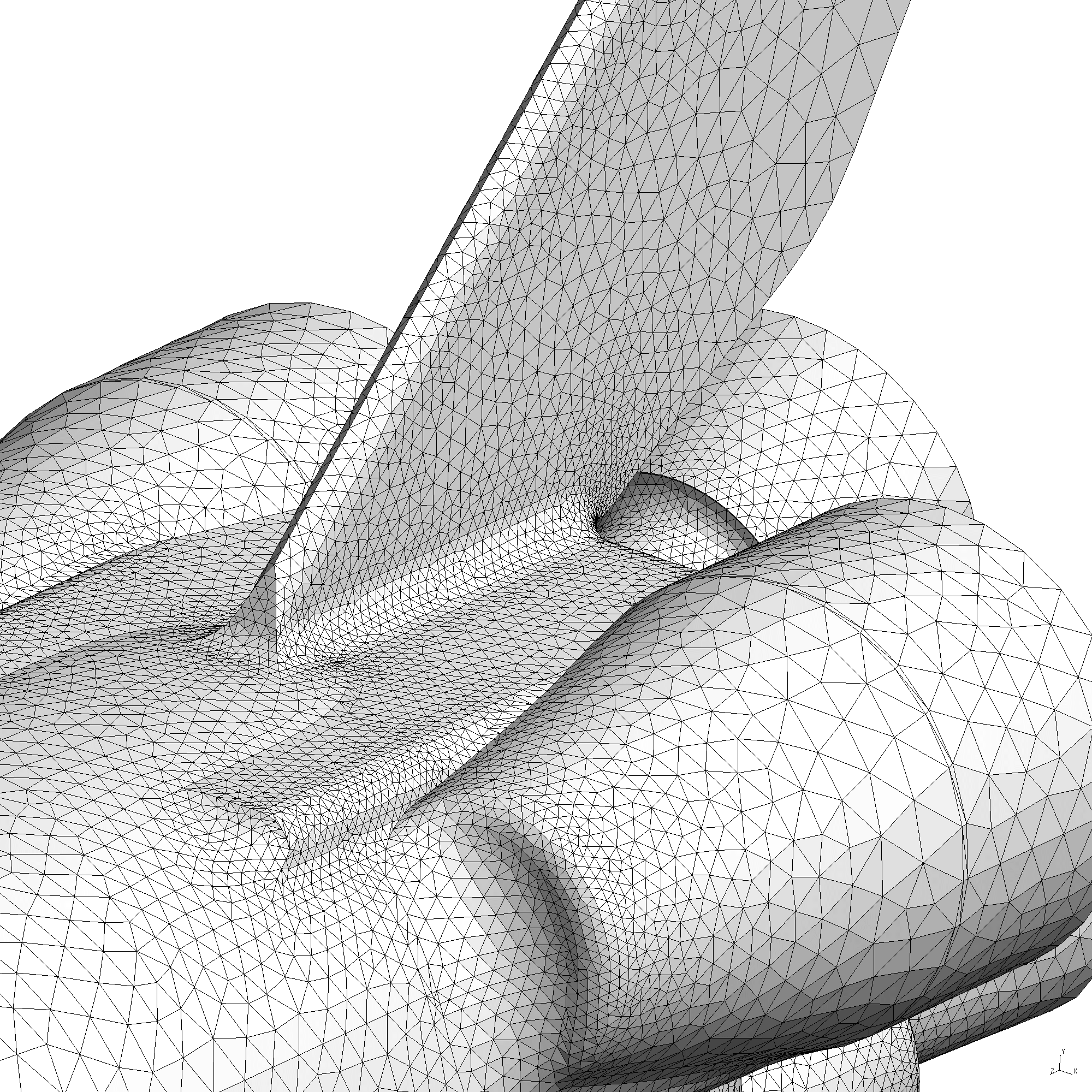}
  \caption{$\alpha = 1.1$}
\end{subfigure}%
\hfill
\begin{subfigure}{.3\textwidth}
  \centering
  \includegraphics[width=\linewidth]{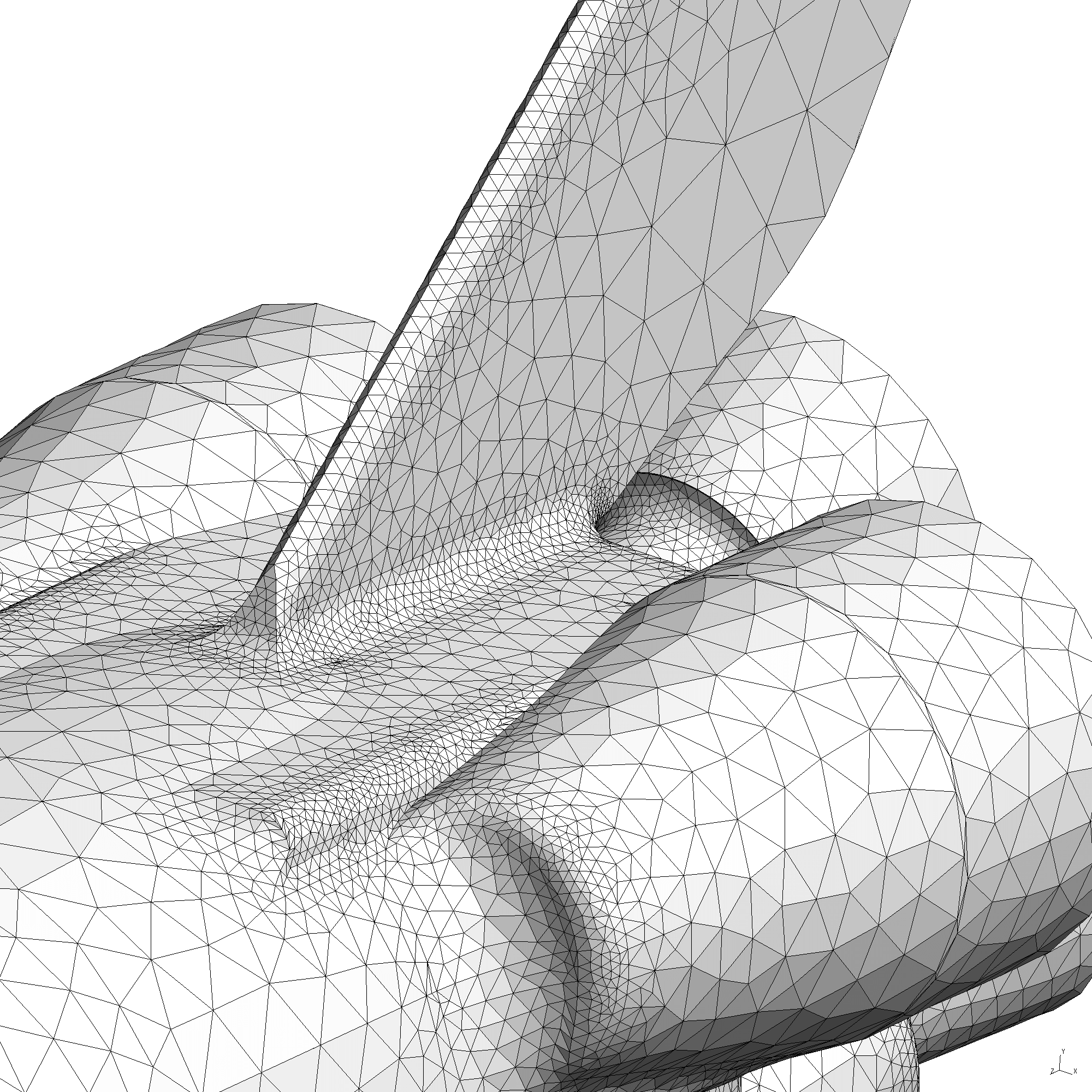}
  \caption{$\alpha=1.4$}
  \label{fig:sub2}
\end{subfigure}
\centering
\hfill
\begin{subfigure}{.3\textwidth}
  \centering
  \includegraphics[width=\linewidth]{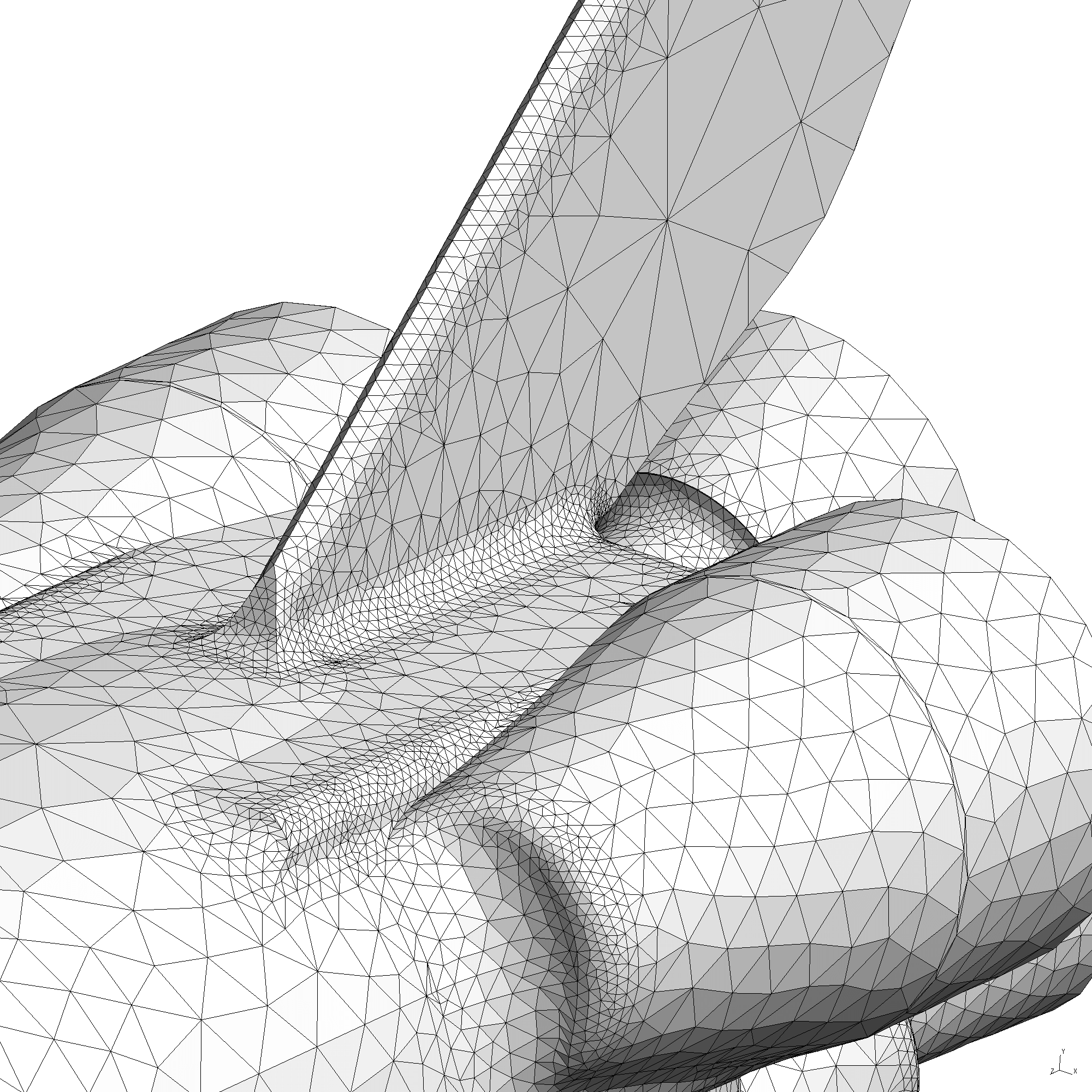}
  \caption{$\alpha = 1.8$}
\end{subfigure}%
\caption{Influence of the gradation $\alpha$ on the final mesh for a node density of $n_d = 20$.}
\label{fig:influenceGradation}
\end{figure}

\emph{Space shuttle.}~This model illustrates the influence of the node density $n_d$ and the gradation $\alpha$ on the resulting mesh. To clearly show the mesh gradation on flat surfaces, surface proximity was not considered to compute these size fields~: only curvature was taken into account. As expected, larger values of $n_d$ will result in more elements on the local osculating circle to a surface, resulting in a finer mesh (Fig.\ref{fig:influenceDensity}). Similarly, a gradation close to $1$ will limit the geometric progression in the mesh size, also resulting in a finer and homogeneous mesh (Fig.\ref{fig:influenceGradation}).
The parameter $n_d$ directly translates into the desired mesh density, and is thus dependent of the aimed application for the mesh. The gradation however, has similar values for different ranges of applications : typically, we recommend using $\alpha$ between $1.05$ and $1.15$ for fluid mechanics computations, $\alpha = 1.2$ for solid mechanics and up to $1.4$ for electromagnetism numerical simulations. The value of $1.4$ is also suited for boundary layer meshes in fluid mechanics. [refs?]

\begin{figure}[h!]
  \centering
  \includegraphics[width=\textwidth,height=\textheight]{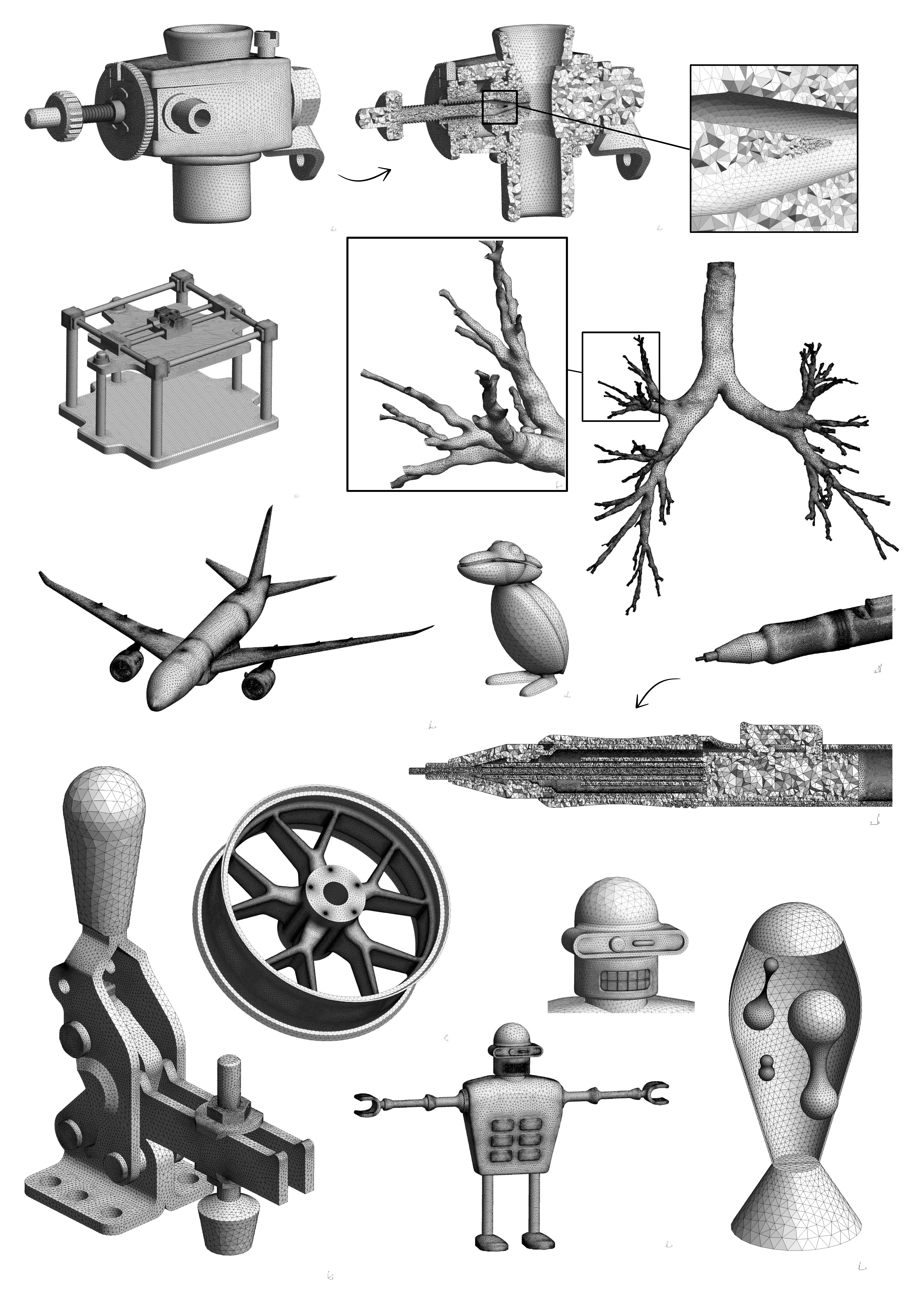}
  \caption{Application of the size field on various geometries from GrabCAD, the ABC Dataset and the Gmsh benchmarks.}
  \label{fig:bunchOfGeometries}
\end{figure}

\emph{Various geometries.}~Finally, we tested our algorithm with the same parameters on various CAD model retrieved from GrabCAD, the ABC Dataset library and the Gmsh benchmarks library (Fig. \ref{fig:bunchOfGeometries}). On all models, accurate size fields were computed in a few seconds to a few minutes, yielding smooth meshes suitable for numerical simulations and saving precious time to the user. Note that feature size computation was not enabled for all models (e.g. for the lava lamp).

\subsection{Adequation between the mesh and the size field}

\begin{figure}
  \centering
  \includegraphics[width=0.49\textwidth]{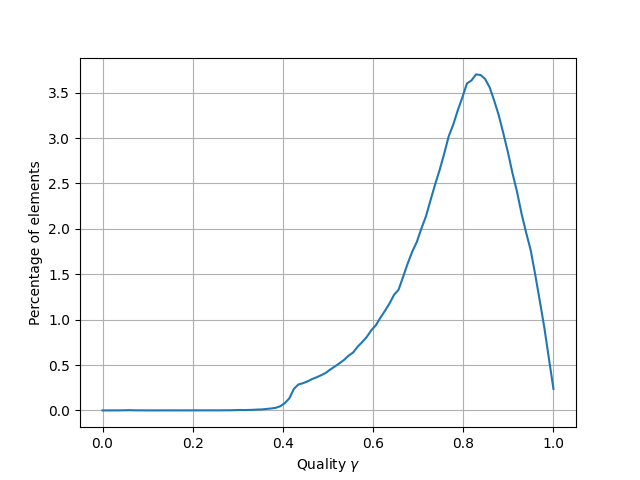}
  \includegraphics[width=0.49\textwidth]{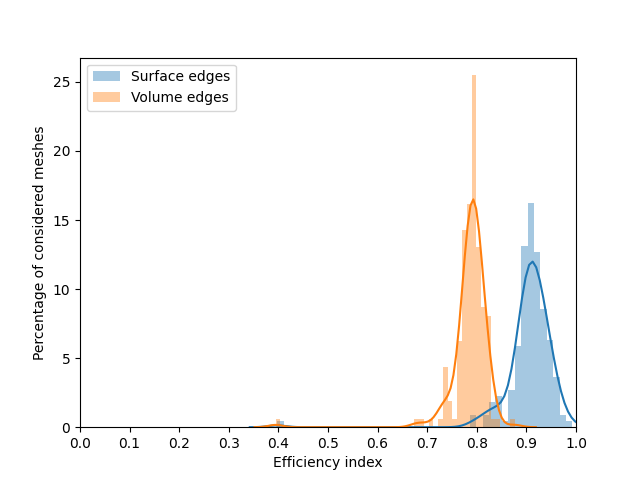}
  \caption{Left : quality of the tetrahedral mesh of the engine block model. Right : efficiency index for a sample of 190 meshes.}
  \label{fig:efficiencyAndQuality}
\end{figure}

\emph{Efficiency index.} ~ The resulting mesh should be as close as possible as a unit mesh in the metric field associated to the isotropic sizing function \cite{frey1999maillages}, i.e. the edges lengths should be as close as 1 when measured in the given metric. Since we are dealing with the discrete counterpart of this continuous mesh framework, we consider that the size specification is respected when the edges lie in $[\frac{1}{\sqrt{2}}, \sqrt{2}]$. To evaluate the adequation between the generated mesh and the prescribed size field, we define an \emph{efficiency index}. Let $l_i$ be the length of an edge $i$ computed in the local metric and $n_e$ the number of edges in the mesh. The efficiency index $\tau$ of a mesh is defined as the exponential to one of all the edges length in the mesh, that is :
$$\tau = \exp{\left(\frac{1}{n_e} \sum_{i=1}^{n_e}\bar{l}_i\right)},$$
with $\bar{l}_i = l_i - 1$ if $l_i < 1$ and $\bar{l}_i = 1/l_i -1$ if $l_i \geq 1$. The efficiency index thus lies in $]0,1]$, the upper bound being reached for a unit mesh. We run our algorithm on close to two hundred CAD models and computed the efficiency index (Fig. \ref{fig:efficiencyAndQuality}), whose median is slightly above $0.8$. This shows that the size field given to the meshing tool is realistic, as it can generate almost unit edges for the given metric field.

\begin{figure}[h!]
\centering
\begin{subfigure}{.5\textwidth}
  \centering
  \includegraphics[width=\linewidth]{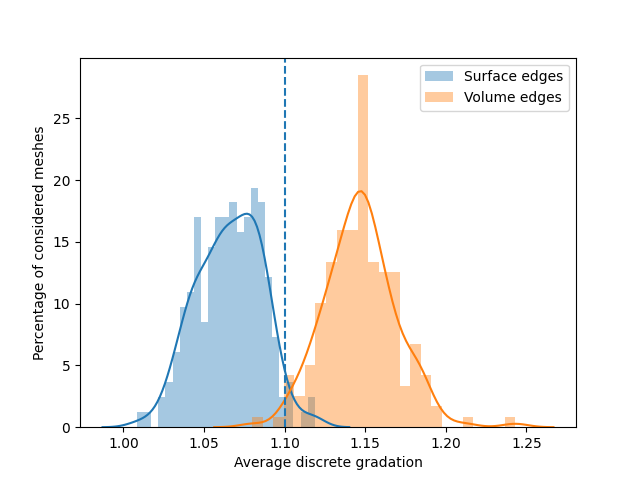}
\end{subfigure}%
\begin{subfigure}{.5\textwidth}
  \centering
  \includegraphics[width=\linewidth]{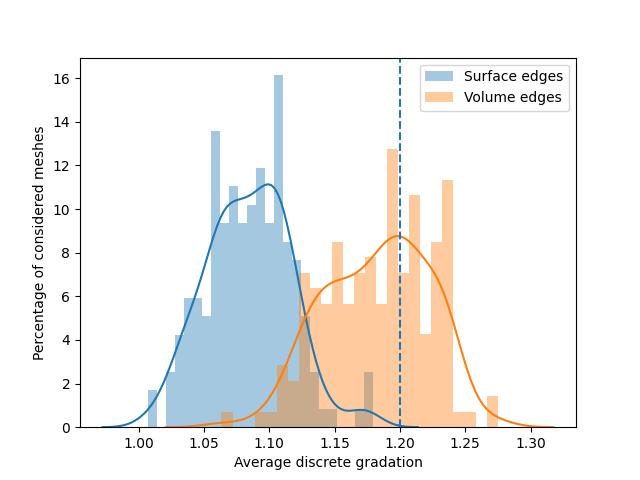}
\end{subfigure}
\caption{Left : discrete gradation for a sample of 190 meshes for a prescribed $\alpha=1.1$. Right : $\alpha=1.2$}
\label{fig:discreteGradation}
\end{figure}

\emph{Discrete gradation.} ~ The resulting mesh should feature a geometric progression in adjacent edges lengths with a ratio up to $\alpha$. Instead of monitoring the maximum edge ratio at every vertices of the final mesh, we use a looser indicator defined on the edges. The \emph{discrete gradation} $\alpha_{d,i}$ measures the progression between the average edges lengths at both vertices of each edge $e_i$. Let $v_1$ and $v_2$ denote the vertices of $e_i$ and $l_{avg}(v)$ the average edge length at vertex $v$. We define the discrete gradation at edge $e_i$ as follows :
\begin{equation*}
  \alpha_{d,i} = \frac{\max [\,l_{avg}(v_1), l_{avg}(v_2)\,] }{\min [\,l_{avg}(v_1), l_{avg}(v_2)\,]}
\end{equation*}
In the same way that the efficiency index is a global indicator for a given mesh, we define an average discrete gradation as the average of $\alpha_{d,i}$ over all edges of the mesh. The average discrete gradation should be close to the user-defined gradation $\alpha$. In our results over the same sample of meshes, the discrete gradation is very close to the required gradation, although it lies slightly above or below depending on the value of $\alpha$ (Fig. \ref{fig:discreteGradation}).

\emph{Quality measure.} ~ Finally, the chosen quality criterion for the tetrahedral meshes is the ratio $\gamma$ of the inscribed radius to the circumscribed radius : $$\gamma = \frac{r}{R}.$$ The resulting mesh for the engine block example contains almost no element with quality below $0.4$, and shows a typical distribution for a 3D mesh.

\subsection{Execution time}
The execution time for the engine block and the Honda engine are given in Table \ref{tab:time}. The size field time is split between the main steps of the algorithm : $(i)$ insert the bounding boxes of the triangles of the surface mesh in the RTree, $(ii)$ compute the approximated curvature with a least square method, $(iii)$ initialize the root octant and refine the octree, $(iv)$ compute the medial axis of the geometry, and $(v)$ limit the size gradient to $\alpha-1$. Computing the medial axis requires the Delaunay tetrahedrization of the points from the surface mesh : this step is detailed and included in the time associated to the medial axis computation. The step "Others" includes all the secondary steps, such as splitting the mesh into faces before computing curvature.

To limit the size gradient, we iterate over the octants until condition \eqref{eq:refinement} is satisfied in all three directions. This step is dependent on the depth of the octree, hence on the minimum size $h_{min}$, the node density $n_d$ and the curvature of the surface mesh given as input. A highly curved region of the model due to a poor resolution of the surface mesh will result in high density in the octree, increasing the smoothing time. This can generally be circumvented by first generating a slightly refined surface mesh, which can be computed quickly and will rule out these extreme curvature magnitudes. Computation time for the medial axis, on the other hand, is linear with the number of nodes of the input mesh : this is indeed observed in Table \ref{tab:time}.

\begin{table}[h!]
 \caption{Execution times}
  \centering
  \begin{tabular}{llrrrr}
    \toprule
     & \multirow{2}[3]{*}{} & \multicolumn{2}{c}{Block} & \multicolumn{2}{c}{Honda engine} \\
    \cmidrule(lr){3-4} \cmidrule(lr){5-6}
     & & Time ($s$) & \% & Time ($s$) & \% \\
    \midrule
    \multirow{8}{*}{Create size field} & Insert surface mesh in RTree & 0.47  & 3.0 \% & 2.84 & 1.3 \%\\
         & Compute curvature          & 0.36  & 2.3  \% & 1.13 & 0.5 \% \\
         & Create and refine octree   & 0.96  & 6.1  \% & 20.59 & 9.1 \% \\
         & Compute medial axis        & 7.19  & 45.7 \% & 39.28 & 17.4 \% \\
         & \emph{incl.} Delaunay tetrahedrization of surface mesh   & 0.61 & 3.9 \% & 3.19 & 1.4 \%\\
         & Limit size gradient  & 5.85  & 37.1 \% & 154.21 & 68.3 \% \\
         & Others               & 0.69  & 4.3  \% & 7.08 & 3.1 \% \\
    \cmidrule(lr){2-6}
         & Total & 15.75 & 100.0 \% & 225.87 & 100.0 \%\\
    \midrule
    \multirow{5}{*}{Mesh model} & Mesh 1D entities & 2.48 & 2.5 \% & 119.80 & 30.2 \%\\
         & Mesh 2D entities & 25.68 & 25.6 \%  & 277.35 & 69.8 \%\\
         & Mesh 3D entities & 72.12 & 71.9 \% & / & / \\
         \cmidrule(lr){2-6}
         & Total & 100.3 & 100.0 \% & 397.33 & 100.0 \% \\
         & \emph{including} Size queries   & 65.13 & 64.9 \% & 169.41 & 42.6 \%\\
    \midrule
    \multirow{5}{*}{} & \# of curves in the CAD model & \multicolumn{2}{c}{1584} & \multicolumn{2}{c}{53,237}\\
    & \# of faces in the CAD model& \multicolumn{2}{c}{533} & \multicolumn{2}{c}{21,771}\\
    \cmidrule(lr){2-6}
    & \# of nodes in input surface mesh & \multicolumn{2}{c}{145,024} & \multicolumn{2}{c}{868,980}\\
    & \# of elements in input surface mesh & \multicolumn{2}{c}{290,116} & \multicolumn{2}{c}{1,737,310}\\
    & \# of octants & \multicolumn{2}{c}{670,566} & \multicolumn{2}{c}{14,615,413}\\
    \cmidrule(lr){2-6}
    & \# of triangles in final surface mesh & \multicolumn{2}{c}{573,808} & \multicolumn{2}{c}{4,605,010}\\
    & \# of tetrahedron in final volume mesh & \multicolumn{2}{c}{2,716,170} & \multicolumn{2}{c}{/}\\
    \bottomrule
  \end{tabular}
  \label{tab:time}
\end{table}

\section{Conclusion}
We have presented a methodology to generate an accurate size field in an automatic way, storing the sizing information in an octree. Five user parameters are required to generate a mesh suitable for numerical simulations, two of which being the minimum and the maximum, or \emph{bulk}, size, whose value is assigned based solely on the characteristic dimension of the CAD model. This leaves the user with only three parameters to choose in order to tune the mesh density. This tool eliminates the tedious operation that is assigning by hand the mesh size on all geometric entities of a CAD model, thus saving precious time. Special care was given to the small features of the geometry : through an approximate medial axis computation, we ensure multiple layers of elements are always generated in narrow regions. This is particularly useful in simulations where Dirichlet boundary conditions are applied. To illustrate our size field computation, we applied our algorithm on a variety of CAD models : in each application, the proposed method accurately targets the regions of higher curvature and adapts the mesh size accordingly. Small features are identified through the medial axis, such that the final meshes include element layers in possible channels in the geometry. For all of the test cases, surface and volume meshes can be obtained in a robust and automatic fashion, and are adapted to the features of the geometric model. This tool will soon be integrated in the standard Gmsh pipeline. Future work will focus on extending to anisotropic size propagation, as well as interaction between the octree and numerical solutions to include \emph{a posteriori} error estimation in the size field design.

\section{Acknowledgements}
The authors would like to thank the Belgian Fund for Scientific Research (FRIA-FNRS) for their support. Financial support from the Simulation-based Engineering Science (Génie Par la Simulation) program funded through the CREATE program from the Natural Sciences and Engineering Research Council of Canada is also gratefully acknowledged.

\bibliographystyle{unsrt}
\bibliography{references.bib}
\end{document}